\documentclass{imsart}

\usepackage{amsthm,amsmath,natbib}

\startlocaldefs

\usepackage{graphics,amssymb,color}

\newcommand{\minimize}{\mathop{\mathrm{minimize}}}

\newcommand{\real}{\mathbb{R}}
\newcommand{\Pp}{{\mathbb P}}
\newcommand{\Ee}{{\mathbb E}}
\newcommand{\lips}{{\cal L}}
\newcommand{\sphere}[1]{O^{#1-1}}
\newcommand{\dimens}{\text{dim}}
\newcommand{\cone}{\text{Cone}}
\newtheorem{theorem}{Theorem}

\def\qed{\hfill $\Box$\newline}

\endlocaldefs

\begin{document}

\begin{frontmatter}

\title{Detecting sparse cone alternatives for Gaussian random fields, with
an application to fMRI}
\runtitle{Detecting sparse cone alternatives}

\author{\fnms{J.E.} \snm{Taylor}\corref{Jonathan E. Taylor}\ead[label=e1]{jonathan.taylor@stanford.edu}\thanksref{t1}}
\thankstext{t1}{Supported in part by NSF grant DMS-0906801.}
\and
\author{\fnms{K.J.} \snm{Worsley}\ead[label=e2]{}\thanksref{t2}}
\thankstext{t2}{Keith Worsley, friend, mentor and colleague passed away February 27, 2009.}

\affiliation{Stanford University, McGill University and University of Chicago}

\address{Department of Statistics\\  Stanford University\\ Sequoia
Hall \\390 Serra Mall\\ Stanford, CA 94305, U.S.A.\\ \printead{e1}}


\runauthor{Taylor \& Worsley}

\begin{abstract}
Our problem is to find a good approximation to the P-value of the maximum of a random field of test statistics for a cone alternative at each point in a sample of Gaussian random fields. These test statistics have been proposed in the neuroscience literature for the analysis of fMRI data allowing for unknown delay in the hemodynamic response. However the null distribution of the maximum of this 3D random field of test statistics, and hence the threshold used to detect brain activation, was unsolved. To find a solution, we approximate the P-value by the expected Euler characteristic (EC) of the excursion set of the test statistic random field. Our main result is the required EC density, derived using the Gaussian Kinematic Formula. 
\end{abstract}



\end{frontmatter}

\section{Introduction}\label{sec:intro}

It seems appropriate to begin this paper with a tribute to the paper's second author, Keith Worsley, for whom this appears posthumously.
This paper is to appear in a volume celebrating David Siegmund's 70th birthday. David and Keith Worsley had worked together several times over their careers \cite{Siegmund:Worsley:1995,Shafie}, most often at the intersection of their two interests: the
distribution of the maximum of random fields. While David's interests
range from the smooth to the non-smooth case, Keith was most interested in smooth random fields
and their application to brain imaging \cite{Worsley:1994:Chi:t:F,friston_statistical_1995,worsley_unified_1996-1}. 
This paper represents  Keith Worsley's last work, before he passed away prematurely
from pancreatic cancer in February 2009. Keith and the first author had discussed this paper right up to a few days before he passed away.  

David has considered two main approaches to such problems: Weyl's volume of tube formulas  as in \cite{johnstone_hotellings_1989,Siegmund:Knowles:1989} and change of measure approaches as in \citet{yakir}. On the other hand, 
Keith preferred using the expected Euler characteristic (EC) approach of \citet{Adler:1981} and his generalizations
\cite{worsley_estimating_1995}. In this paper, we combine the EC approach to the volume of tube formula
via the Gaussian Kinematic Formula \citep{Taylor:2002:Kinematic} which expresses Keith's EC densities in terms of coefficients the {\em Gaussian measure} of a tube. Referring back to David Siegmund's approach to these problems, these coefficients are also coefficients
in an expansion of their own  change of measure formula on Gaussian space \citep{taylor_vadlamani}. 

This paper is concerned with the maxima of (functions of) smooth Gaussian random fields. Let $T(s)$, $s\in\real^D$ be a random field, and let $S\subset\real^D$ be a fixed {\em search region}. Our main interest is to find good approximations to the P-value of the maximum of $T(s)$ in $S$:
\begin{equation}\label{eq:Pmax}
\mathbb{P}\left(\max_{s\in S} T(s)\geq t\right).
\end{equation}
The random field $T(s)$ will be one of a variety of test statistics for a cone alternative in a multivariate Gaussian random field. Two of these test statistics have been proposed in the neuroscience literature \citep{Friman:2003,Calhoun:2004} but without a P-value (\ref{eq:Pmax}). \citet{Worsley:Taylor:2005:latency} gives a heuristic approximation to the P-value of the \citet{Friman:2003} statistic. This has been incorporated into the {\tt R} package {\tt fMRI} \citep{fmriR}. This paper aims to give a correct P-value approximation to both of these test statistics and the likelihood ratio test statistic.

To do this, we first define the test statistic random fields in Section \ref{sec:test}, then evaluate their approximate P-values (\ref{eq:Pmax}) in Section \ref{sec:pvalue} using the EC heuristic and the Gaussian kinematic formula. Section \ref{sec:pvalue} concludes with a subsection that relates our methods to those we have used for the Hotelling's $T^2$ random field \citep{Taylor:Worsley:Roy}. Finally in Section \ref{sec:app} we apply our methods to the re-analysis of an fMRI data set already used for the same purpose in \citet{Worsley:Taylor:2005:latency}.

\section{The test statistics}\label{sec:test}

\subsection{Definitions of the test statistics}

The test statistics are defined as follows.
Let $Z(s)=(Z_1(s),\dots,Z_n(s))'$, $s\in S\subset \real^D$, be a vector of $n$ i.i.d. Gaussian random fields with
$$
\mathbb{E}(Z(s))=\mu(s),\ \ \ \ \mathbb{V}(Z(s))=\sigma(s)^2I_{n\times n}.
$$
Usually $\sigma(s)$ is unknown and must be estimated separately at each point, so keeping this in mind, we will set $\sigma(s)=1$ without loss of generality.
Let $U\subset \sphere{n}$, the unit $(n-1)$-sphere. At each $s \in S$,
we are interested in testing that the mean is zero against the cone alternative:
$$
H_{0}: \mu(s)=0\ \ \ \text{vs.}\ \ \ H_{1}:\mu(s) \in \cone(U) =
\{c \cdot u: c \geq  0,\, u \in U\}
$$
\citep{Rob:Wri:Dyk}. The
likelihood ratio test of $H_{0}$ vs. $H_{1}$ is equivalent to
\begin{equation}
  \label{eq:bar:chi}
\bar\chi(s)=\max_{u\in U} u'Z(s),
\end{equation}
which we call the {\em $\bar\chi$ random field} because it has a
so-called $\bar\chi$ marginal
distribution when $\cone(U)$ is convex (see Section
\ref{sec:bar:mix} below).
As mentioned above, $\sigma(s)$ is usually unknown so the $\bar\chi$ random field must be normalized separately at every point $s$. We shall consider two ways of doing this.

The first is the {\em likelihood ratio cone random field}, equivalent to the likelihood ratio of the cone alternative under unknown variance:
$$
T_{\rm LR}(s)=\frac{\bar\chi(s)}{\sqrt{(||Z(s)||^2-\bar\chi(s)^2)/n}},
$$
or equivalently, the maximum correlation between a point in the cone
and the data.
The second, proposed by \citet{Friman:2003}, is only defined if
$U$ is a subset of some $k$-dimensional subspace of $\real^n$,
in which case there are effectively $\nu=n-k$ residual degrees of freedom
 which can be used to estimate $\sigma(s)$ and
normalize
$\bar\chi(s)$.
 Suppose $Z_\bot(s)$ is the projection of $Z(s)$ onto the orthogonal
 complement of the linear span of $U$, so that $Z_\bot(s)$ is
 independent of $\bar\chi(s)$ and has mean 0 under $H_1$. Then the {\em independently normalized cone random field} is
$$
T_{\rm IN}(s)=\frac{\bar\chi(s)}{||Z_\bot(s)||/\sqrt{\nu}}.
$$
Note that if $U=\sphere{k}$ (by this we mean a $(k-1)$-sphere embedded in $\real^n$) then the two cone random fields are both equivalent to the {\em F-statistic random field}
$$
F(s)=\frac{||Z_\top(s)||^2/k}{||Z_\bot(s)||^2/\nu}
$$
where $Z_\top(s)$ is the projection of $Z(s)$ onto the linear subspace
spanned by $U$.

%

\begin{figure}[t]\begin{center}
\mbox{\scalebox{1}{\includegraphics*{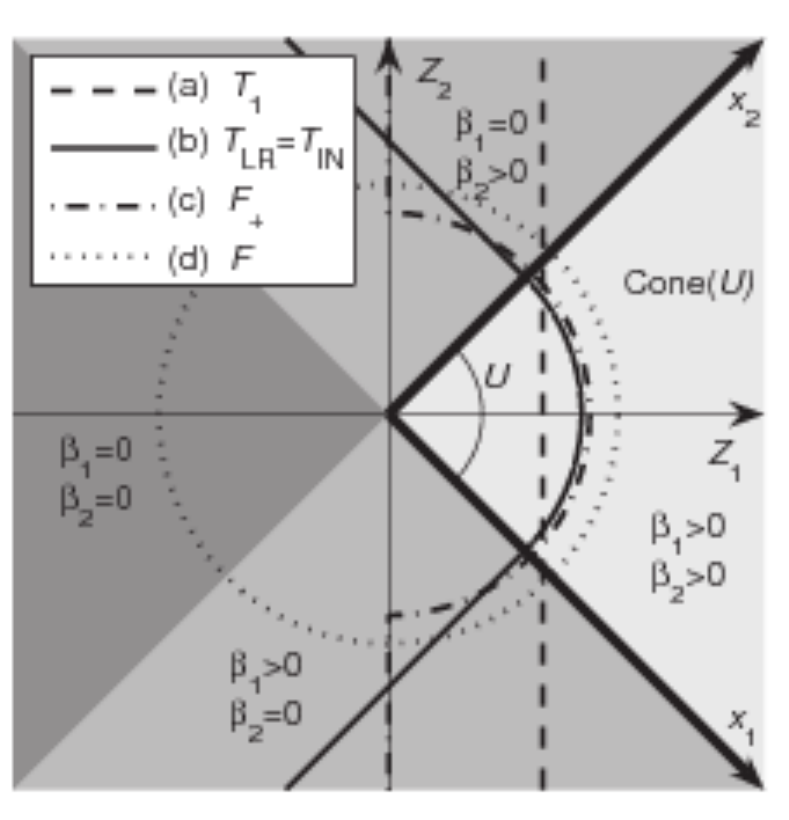}}} \end{center}
\caption{Rejection regions
(the side of the boundary that excludes the origin) of the test statistics at
$P=0.05$ with infinite sample size for a 2D ($k=2$) right-angled cone alternative covering the first two components $Z_1, Z_2$ of $Z$. The middle of the cone $u$ is parallel to the $Z_1$ axis. The cone can also be expressed as a linear model with $m=2$ regressors $x_1$ and $x_2$ with non-negative coefficients $\beta_1\geq 0$ and $\beta_2\geq 0$. The $\bar\chi$ statistic is the length of the projection of $Z$ onto the nearest edge of the cone (including the vertex of the cone and the interior of the cone itself). The null distribution of $\bar\chi$ is a mixture of $\chi_j$ random variables with weights $p_j=\Pp(\#\{\widehat\beta{\rm 's}\geq 0\}=j)$ equal to the relative size of the shaded regions: $p_{0,1,2}=1/4, 1/2,1/4$.
The statistic $F_+$ is the one-sided $F$ statistic of \citet{Calhoun:2004}.}\label{figrej}
\end{figure}

\subsection{Power and maximum likelihood}

Both cone statistic random fields should be more powerful than the
F-statistic random field since the F-statistic wastes power on
alternatives that are outside the cone. The one-sided F-statistic
tries to make up for this, but it is inadmissable (for
infinite $\nu$ and fixed $s$) because its acceptance region is concave
\citep{Birnbaum:1954} - see Figure \ref{figrej} - although it is not
clear how to construct a test which
dominates it. If in fact the alternative is at the middle
of the cone then $T_1$ should be the most powerful.

Between the two
cone statistics, the advantage of $T_{\rm LR}(s)$ is that it uses all
the information in the data to estimate the variance and so it should
be more powerful than $T_{\rm IN}(s)$.
\citet{Cohen:Sack:Inadmiss} show that $T_{\rm LR}(s)$ is admissible in specific examples, whereas $T_{\rm IN}(s)$ is always inadmissable.
However if in fact the mean is
outside the cone but still inside the linear subspace spanned by $U$,
then we would expect $T_{\rm IN}(s)$ to be more powerful. The reason
is that a mean $\mu(s)$ outside the cone would increase the
denominator of $T_{\rm LR}(s)$ but not that of $T_{\rm
  IN}(s)$. \citet{Friman:2003} chose the more conservative $T_{\rm IN}(s)$. This strategy sacrifices a few degrees of freedom and a small loss of power if $\mu(s)$ really is in the cone, against a much larger loss of power if it is not. \citet{Worsley:Taylor:2005:latency} investigates power in an fMRI application that we shall also use in Section \ref{sec:app}. For a general discussion of power and likelihood ratio tests in this setting see \citet{Emperor:New:Tests}.


We note in passing that we have used maximum likelihood principles only at a
single point $s$, not over the whole space $S$, which would require a spatial
model for the mean and covariance function of the random fields. In the case of
known $\sigma(s)$, a standard reproducing kernel argument, discussed  in \citet{Siegmund:Worsley:1995},
can be used to show that if each of the components of $\mu(s)$ is proportional to the
spatial correlation function centered at some unknown point $s_0$
(which is assumed to be the
same for each component),  then $\max_{s \in S}\bar\chi(s)$ is the likelihood ratio test statistic.


Our interest is confined to $s$ in a search region $S\subset\real^D$, where we expect $H_0$ to be true at most points, with only a sparse set of points $S_1$ where $H_1$ is true.
This suggests that we should estimate $S_1$ by thresholding the above test statistic
random fields at some suitably high threshold. Choosing the threshold which controls the P-value of the maximum of
the random field to say $\alpha=0.05$ should be powerful at detecting $S_1$, while
controlling the false positive rate outside $S_1$ to something
slightly smaller than $\alpha$. Our main problem is therefore to
find the P-value of the maximum of these random fields of test statistics (\ref{eq:Pmax}), which is the main aim of this paper.

\subsection{Mixture representation of $\bar\chi$}
\label{sec:bar:mix}

The $\bar\chi$ random field is so-named because it has a useful representation
in terms of a mixture of $\chi_j$ random fields with $j$ degrees of freedom \citep{Lin:Lindsay:1997,Takemura:Kuriki:Chi:Bar}.
The mixture representation works when $\cone(U)$ is convex and
polyhedral, and asymptotically when $\cone(U)$ is only locally convex (see
Section \ref{sec:bar} below). The simplest way of seeing where the
polyhedral cone enters the picture is to write it as a linear model with non-negative coefficients:
\begin{equation} \label{eq:linmod}
H_1: \mu(s) =\sum_{j=1}^m x_j\beta_j(s), \ \ \beta_{1}(s),\dots,\beta_m(s)\in\real^+.
\end{equation}
The regressors $x_1,\dots,x_m\in \real^n$ contain the vertices of $U$ (times arbitrary scalars), and they may be linearly dependent (see Figure \ref{figrej}). The cone may even contain linear subspaces  (for instance, take
$x_2=-x_1$ above) which effectively corresponds to having a certain
number of unrestricted coefficients in $\mu(s)$ under $H_1$.

To actually compute the $\bar\chi(s)$ random field, one must solve a convex
problem at each location $s$. This can be done in several ways: the most direct is to first perform
all-subsets least-squares regression, then throw out any fitted
model that has negative coefficients. Amongst those
that are left, the model that fits the best, with fitted values
\begin{equation} \label{eq:fitlinmod}
\widehat Z(s)=\widehat \mu(s)=\sum_{j=1}^m x_j\widehat \beta_j(s), \ \ \widehat\beta_{1}(s),\dots,\widehat\beta_m(s)\in\real^+,
\end{equation}
is the maximum likelihood
estimator of $\mu(s)$, and $\bar\chi(s)=||\widehat Z(s)||$. Alternatively, one may solve the problem
\begin{equation}
\minimize_{(\beta(s))_{s \in S}} \sum_{s \in S} \|Z-X\beta(s)\|^2_2 \ \ \text{subject to} \ \ \beta_i(s) \geq 0, \ 1 \leq i \leq m, s \in S.
\end{equation}
This is is a collection of separable convex problems, each of which can be solved via coordinate descent \cite{friedman_pathwise_2007} or first-order methods (c.f. \cite{nesta}). As the inputs are smooth, one would expect that
warm starts at adjacent locations would greatly speed up the convergence of such algorithms.
There is a huge literature on such non-negative least squares (NNLS) problems, with many applications in inverse problems, and many faster algorithms than all-subsets regression, such as the classic one by \citet{Lawson:Hanson}.

From a geometric perspective, estimation of $\mu(s)$ is
equivalent to projecting $Z(s)$ onto $\cone(U)$, i.e., finding the
face of $\cone(U)$ closest to $Z(s)$. Here, a face of $\cone(U)$ could represent
the vertex of $\cone(U)$, in which case $\widehat Z(s)=0$; an edge of
$\cone(U)$; or even the interior of $\cone(U)$, in which case
$\widehat Z(s)=Z(s)$.
Let $A\subset \cone(U)$ represent a
generic face of $\cone(U)$.
Further, let $\widehat Z_A(s)$ be the projection of $Z(s)$ onto the
linear subspace spanned by $A$, so that $\{\widehat Z_A(s)\in
\cone(U)\}$ is the event that the non-negativity restrictions are satisfied for face
$A$. Then,
\begin{equation} \label{eq:bar:chi:pre}
\bar\chi(s) = \max_{A}1_{\{\widehat Z_A(s)\in
\cone(U)\}}\cdot \|\widehat Z_A(s)\|,
\end{equation}
and let $\widehat A(s)$ be the value of $A$ that achieves this
maximum. Actually, there are values of $Z(s)$ for which more than one face achieves
the maximum above, though these occur on lower dimensional subsets of
$\real^n$, which correspond to lower dimensional surfaces in the
search region $S$.
From \eqref{eq:bar:chi:pre}, it is clear that
\begin{equation} \label{eq:bar:chi:dec}
\bar\chi(s) = \sum_{A}1_{\{\widehat A(s)=A\}} \cdot \|\widehat Z_A(s)\|.
\end{equation}
Clearly,
$$
\bar\chi(s) \bigl|\{\widehat A(s)=A\} \sim \chi_{\dimens(A)},
$$
which only depends on the dimensionality of $A$, and so
$$
\bar\chi(s) \bigl|\{\dimens(\widehat A(s))=j\} \sim \chi_j.
$$
Hence its unconditional marginal distribution is a mixture of $\chi_j$'s
\begin{equation} \label{eq:Pchibar}
\Pp(\bar\chi(s) \geq t) = \sum_{j=0}^n p_j(U) \Pp(\chi_j \geq t)
\end{equation}
with weights
$$
p_j(U)=\Pp\left(\dimens(\widehat A(s))=j\right),\ \ \ \ \ 0 \leq j \leq n.
$$
These weights are the probability that the face of $\cone(U)$ that is closest to
$Z$ has dimension $j$, or, in terms of the fitted linear model (\ref{eq:fitlinmod}),
$$
p_j(U)=\Pp\left(\#\{\widehat \beta{\rm 's}>0\}=j\right),\ \ \ \ \ 0 \leq j \leq n.
$$
Above,
we define $\chi_0=0$ to be a constant random variable
which corresponds to $Z(s)$ being closest to the vertex of $\cone(U)$.
Depending on the structure of $\cone(U)$, one or more of the
$p_j(U)$'s may be zero. More specifically, let $L(U)$ be the largest linear subspace contained in
$\cone(U)$ with $L(U)$ possibly equal to 0, the  subspace containing
only the 0 vector.
It is not hard to see that
$$
 l(U) \overset{\Delta}{=} \dimens(L(U)) = \min\{j: p_j(U) > 0\}$$
and further,
$$ \|\widehat Z_{L(U)}(s)\| \leq \bar\chi(s) \leq \|Z(s)\|.$$
Finally, we also note that, for $t > 0$, $\Pp(\chi_0 \geq  t)=0$ so effectively the sum in
\eqref{eq:Pchibar} is really a sum over $1 \leq j \leq n$ and
we can generally ignore $p_0(U)$ which we do in later expressions for
the EC densities of $T_{\rm IN}(s)$ and $T_{\rm LR}(s)$.

By
approximation, this argument extends to general convex cones, though
the $p_j$'s have slightly different interpretations even though they
are limits of the $p_j$'s of the polyhedral approximations, see
Section \ref{sec:bar} below \citep{Lin:Lindsay:1997,Takemura:Kuriki:Chi:Bar}.

\begin{figure} \begin{center}
\mbox{\scalebox{0.8}{\includegraphics*{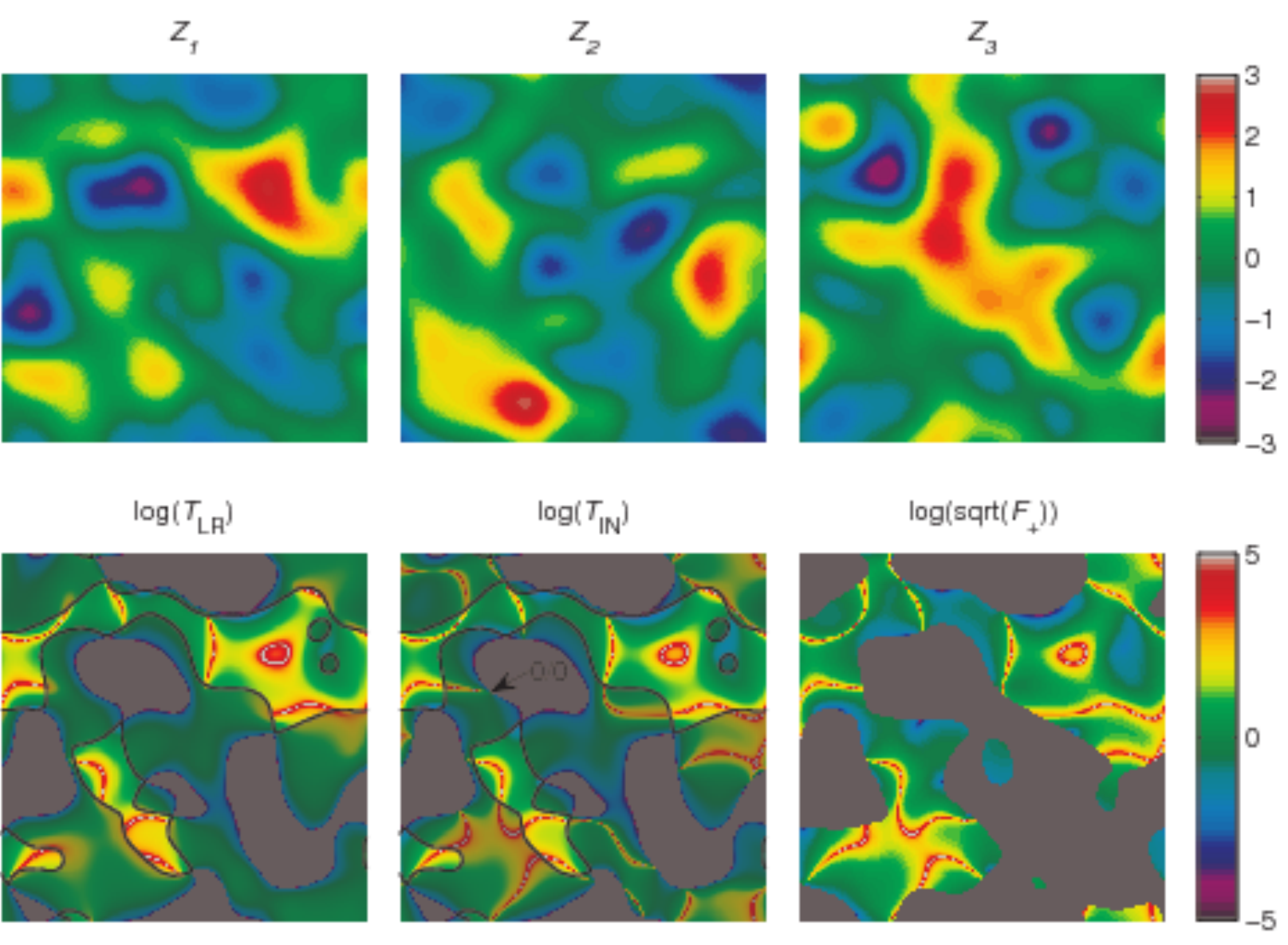}}}\end{center}
\caption{Examples of $n=3$ Gaussian random fields in $D=2$ dimensions (top row). Bottom row: the random fields $T_{\rm LR}$, $T_{\rm IN}$ and $F_+$ for the same quarter circle cone as in Figure \ref{figrej}, so that $k=2$ and $\nu=1$.
In the three patches the $\bar\chi$ random fields are $\chi_j$ fields with $j=$ dimensionality of the nearest cone face. In the gray patches, $j=0$, $T_{\rm LR}=T_{\rm IN}=F_+=0$; in the medium shaded patches, $j=1$, $T_{\rm LR}^2\sim F_{1,2}$ and $T_{\rm IN}^2\sim F_{1,1}$; in the unshaded patches, $j=2$, $T_{\rm LR}^2=T_{\rm IN}^2=F_+\sim F_{2,1}$ (times scalars). The boundary between the medium shaded and unshaded patches (heavy black line) is the edge of the cone, $x_1$ or $x_2$. When the denominator has one degree of freedom, the statistic takes the value $\infty$ on random curves; when it has two degrees of freedom, it takes the value $\infty$ only at the points where these curves touch the boundary. $T_{\rm IN}$ is not defined everywhere because it takes the value 0/0 at random points (arrow). }\label{fig0d0}
\end{figure}

Note that while the marginal distribution of the $\bar\chi(s)$ random
field is a mixture
of $\chi_j$ random variables, it is not strictly a mixture as a random
field. Rather, realizations of the random field resemble a {\em patchwork}
of $\chi_j$ random fields with patches $\{s: \widehat{A}(s)=A\}$ on
which
we observe $\|\widehat{Z}_A(s)\|\sim\chi_{{\rm dim}(A)}$  (see Figure \ref{fig0d0}).

This representation also sheds some light on the
two normalized random fields $T_{\rm LR}(s)$ and $T_{\rm IN}(s)$ as
patchwork mixtures of $\sqrt{F}$ random fields of appropriate degrees of freedom.
In terms of the representation \eqref{eq:bar:chi:dec},
it is not hard to see that
\begin{equation}
  \label{eq:LR:dec}
  T_{\rm LR}(s)= \sum_{A}1_{\{\widehat A(s)=A\}}
  \cdot \frac{\|\widehat Z_A(s)\|}{\|Z(s)- \widehat Z_A(s)\|/\sqrt{n}}.
\end{equation}
Above, some slight care must be taken at points $s$ contained in the intersection of the
closure
of two or more patches. For these points, we can arbitrarily assign
$\widehat{A}(s)$
to any appropriate face of $\cone(U)$.
The representation \eqref{eq:LR:dec} shows immediately that its marginal distribution is that of a
mixture of $\sqrt{jn/(n-j) \cdot F_{j,n-j}}$ random variables with
weights $p_j(U)$. As in the $\chi_0$ case, we define $F_{0,l}=0$ to be
a constant random variable for all $l$. For the independently normalized cone random field
\begin{equation}
  \label{eq:IN:dec}
\begin{aligned}
  T_{\rm IN}(s)= \sum_{A}1_{\{\widehat A(s)=A\}}
  \cdot \frac{\|\widehat Z_A(s)\|}{\|Z_{\bot}(s)\|/\sqrt{\nu}}
\end{aligned}
\end{equation}
which shows that its marginal distribution is a mixture of $\sqrt{j
\cdot F_{j,\nu}}$ random variables with weights $p_j(U)$.

\subsection{Dimensionality}\label{sec:dim}

The representation of $T_{\rm IN}(s)$ and $T_{\rm LR}(s)$ as
patchwork mixtures of $\sqrt{F}$
random fields shows that we must consider constraints on $D$ dictated by
the total degrees of freedom $n$ and $\cone(U)$ (see Figure \ref{fig0d0}). For the $F$
random field, recalling the argument in \citet{Worsley:1994:Chi:t:F}, we note that the
set where $||Z(s)||$ takes the value zero is the intersection of the
zero sets of each of the components of $Z(s)$, so its dimensionality is $D-n$
if $D\geq n$ or empty if $D<n$. This means that if $D\geq n$ then
$F(s)=0/0$ with positive probability somewhere inside $S$, in which
case $F(s)$ is not defined. Hence we must have $D<n$ for $F(s)$ to be
well defined. The same argument applies to
$F_+(s)$ and to $T_1(s)$ for which we must have $D < \nu+1$.

By a similar argument, $T_{\rm LR}(s)$ is made up of $\sqrt{F_{j,n-j}}$
random fields for $l(U) \leq j \leq n$, so we must have $D < n$ to avoid
0/0 for such random fields. A similar argument applies to $T_{\rm IN}(s)$
though the limit on the dimension is more restrictive and slightly
more difficult to describe. In principle, we simply want to avoid 0/0
for the random field $T_{\rm IN}(s)$. However, when $l(U)=0$, we can allow some isolated
0/0 points
 within the interior of the patch $\{s: \widehat{A}(s)=0\}$, i.e. when the numerator of
 $T_{\rm IN}(s)$
is 0. If we allow more than isolated points, say curves of 0/0, these
will generally intersect the boundary of the patch
$\{s:\widehat{A}(s)=0\}$
causing $T_{\rm IN}(s)$ to be undefined at such points (see the white arrows in Figure \ref{fig0d0}(a,b)).
In other words, we really need to avoid 0/0 on the closure of the set $\{s: \widehat{A}(s) \neq
0\}$. When $l(U)=0$, on this set
$$
\min\{\|\widehat{Z}_A(s)\|: \dimens(A)=1\} \leq \bar\chi(s) \leq
\|Z(s)\|$$ therefore there will be no 0/0's if there are no 0/0's for any of
the $F_{1,\nu}$ random fields
$$
\left\{\frac{\|\widehat{Z}_A(s)\|^2}{\|Z_{\bot}(s)\|^2/\nu}: \dimens(A)=1\right\},$$
that is, if $D < \nu+1$. However, if $l(U) > 0$, then
$\{s:\widehat{A}(s) = 0\}$ is of strictly lower dimension than $D$
and even isolated 0/0 points within this patch will cause $T_{\rm
  IN}(s)$ to be undefined, hence we must again avoid 0/0's in the
closure of $\{s: \widehat{A}(s) \neq 0\}$ which is just $S$, the
entire search region.
As noted in the previous section, when $l(U) > 0$
$$
\|\widehat{Z}_{L(U)}(s)\| \leq \bar\chi(s) \leq \|Z(s)\|$$
and there will be no 0/0's in $T_{\rm IN}(s)$ if there are no 0/0's in
the $F_{l(U),\nu}$ random field
$$
\frac{\|\widehat{Z}_{L(U)}(s)\|^2/l(U)}{\|Z_{\bot}(s)\|^2/\nu},$$
that is, if $D<\nu+l(U)$. In summary, considering both cases
$l(U)=0$ and $l(U)>0$, we must have
$D<\nu+\max(l(U),1)$.

When $\cone(U)$ is non-convex, the situation is more difficult to
describe in exact terms for both $T_{\rm IN}(s)$ and $T_{\rm LR}(s)$. If $\cone(U)$ is non-convex, then the marginal distribution of
$\bar\chi(s)$ is no longer exactly a mixture of $\chi_j$'s with the
error being exponentially small \cite{TTA:Validity}. 

\section{P-value of the maximum of a random field}
\label{sec:pvalue}

\begin{figure}[t] \begin{center}
\mbox{\scalebox{0.75}{\includegraphics*{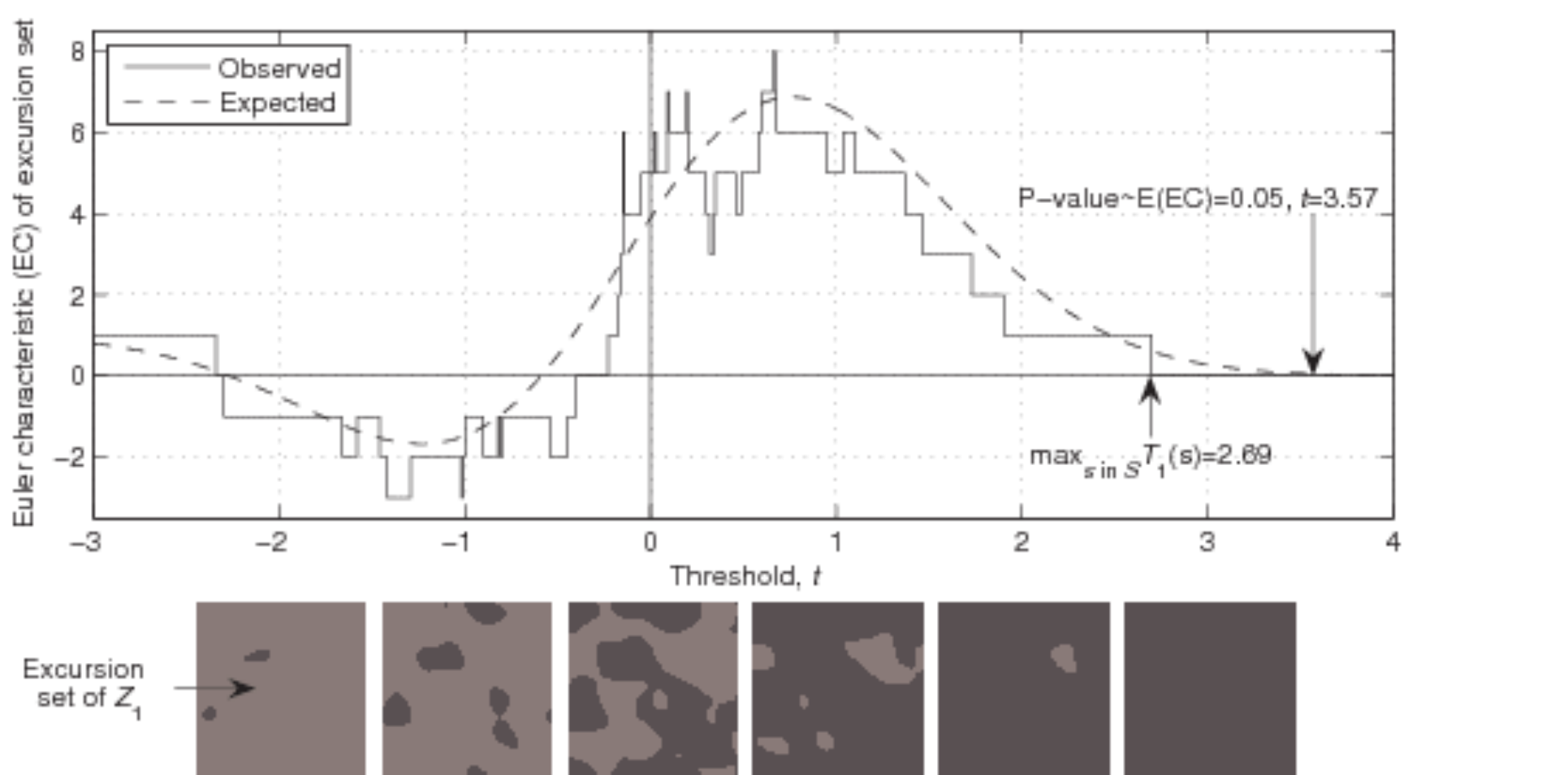}}}
\end{center}
\caption{The Euler characteristic (EC) of excursion sets of the Gaussian random field $Z_1$ from Figure \ref{fig0d0} plotted against threshold $t$, together with the expected EC under $H_0$ from (\ref{p}). Bottom row: the excursion sets (light gray) for $t=-2,\dots,3$; the search region $S$ is the whole image. At high thresholds the expected EC is a good approximation to the P-value of the maximum (arrowed). The approximate $P=0.05$ threshold is $t=3.57$ (arrowed).}\label{figEC}
\end{figure}

A very accurate approximation to the P-value of the maximum of any smooth
isotropic random field $T(s)$, $s\in S \subset \real^D$, at high thresholds
$t$, is the expected Euler characteristic (EC) $\varphi$ of the excursion set:
\begin{equation}\label{p}
{\mathbb P}\left(\max_{s\in S} T(s) \geq t \right) \approx {\mathbb
E}(\varphi\{s\in S : T(s) \geq t \}) = \sum_{d=0}^D \lips_d(S)\rho_d(t),
\end{equation}
where $\lips_d(S)$ is the $d$-dimensional {\em intrinsic volume} of $S$
(defined in Appendix \ref{app:iv}), and $\rho_d(t)$ is the $d$-dimensional {\em
EC density} of the random field above $t$
\citep{Adler:1981,Worsley:1995:Boundary,Adler:2000,RFG}. The heuristic is that for high thresholds the EC takes the value 0 or 1 if the excursion set is empty or not, so that the expected EC approximates the P-value of the maximum (see Figure \ref{figEC}). The approximation is extraordinarily accurate, giving exponential accuracy for Gaussian random fields \citep{TTA:Validity}.
A different approach using volumes of tubes \citep{Siegmund:Knowles:1989,Johansen:Johnstone:1990,Sun:1993,Sun:Loader:1994,Sun:2000,Pilla:2006} is, in our context, essentially the same as the methods used here, as shown by \citet{Takemura:Kuriki:Equivalence:2002}.

For $D=3$, our main interest in applications, $\lips_{0,1,2,3}(S)$ are: the EC,
twice the `caliper diameter', half the surface area, and the volume of $S$
respectively (for a convex set, the caliper diameter is the average distance between the two parallel tangent planes to the set). If the random field $T(s)$ is a function of Gaussian random fields, such as all the test statistic random fields considered so far, and these Gaussian random
fields are non-isotropic, then it is only necessary to replace intrinsic volume
in (\ref{p}) by {\em Lipschitz-Killing curvature}. Lipschitz-Killing curvature
depends on the local spatial correlation of the component Gaussian random
fields, as well as the search region $S$
\citep{Taylor:Adler:2003:EC:Manifolds,TW:noniso}.

Morse theory can be used to obtain the EC density of a smooth random field $T=T(s)$ as
\begin{equation}\label{rho}
\rho_d(t)={\mathbb E}\left(1_{\{T\geq t\}}{\rm det}(- \ddot T_d)\ |\ \dot
T_d=0\right){\mathbb P}(\dot T_d = 0),
\end{equation}
where dot notation with subscript $d$ denotes differentiation with respect to
the first $d$ components of $s$ \citep{Worsley:1995:Boundary}. For $d=0$, $\rho_0(t)={\mathbb P}(T\geq t)$.
The Morse method of obtaining EC densities, though straightforward in principle, usually involves an enormous amount of tedious algebra. Entire papers have been devoted to evaluating (\ref{rho}) for an ever wider class of random fields of test statistics such as Gaussian \citep{Adler:1981}, $\chi^2$, $T$, $F$ \citep{Worsley:1994:Chi:t:F}, Hotelling's $T^2$ \citep{Cao:Worsley:1999:Hotelling}, correlation
\citet{Cao:Worsley:1999:Correlation}, scale space \citep{Siegmund:Worsley:1995,Worsley:Scale:Chisq:2001,Shafie} and Wilks's $\Lambda$ \citep{Carbonell:Wilks}. A much simpler method is given in the next section.

\subsection{The Gaussian Kinematic Formula}

There is a much simpler way of getting EC densities when $T$ is built from
independent {\em unit Gaussian random fields} (UGRF). A UGRF is a Gaussian random field with zero mean,
unit variance, and identity variance of its spatial derivative. Note that any
stationary Gaussian random field can be transformed to a UGRF by appropriate
linear transformations of its domain and range. Without loss of generality we shall assume that all the random fields considered so far are built from UGRFs.

This simpler method is based on the {\em Gaussian Kinematic Formula} discovered by \citet{Taylor:2002:Kinematic}. The idea is to take
the Steiner-Weyl volume of tubes formula (\ref{eq:voltub}) and replace the
search region by the rejection region, and volume by probability. Somewhat miraculously, the
coefficients of powers of the tube radius are (to within a constant) the EC
densities we seek.

The details are as follows. Suppose $T(s)=f(Z(s))$ is a function of UGRFs $Z(s)=(Z_1(s),\dots,Z_n(s))'$. Put a tube of radius $r$ about the rejection
region $R_t=\{z \in \real^n: f(z)\geq t\}\subset \real^n$, evaluate the
probability content of the tube (using the ${\rm N}_n(0,I_{n\times n})$
distribution of $Z=Z(s)$), and expand as a formal power series in $r$. Denoting the tube by
${\rm Tube}(R_t,r)=\{x: \min_{z\in R_t}||z-x||\leq r\}$, then
\begin{equation}\label{gkf}
\Pp\left( Z \in {\rm Tube}(R_t,r)\right)=\sum_{d=0}^\infty \frac{r^d}{d!}
(2\pi)^{d/2}\rho_d(t).
\end{equation}
Since the spatial dependence on $s$ is no longer needed, we omit it until further notice.

For example, let $f(z) = u'z$ for fixed $u$ with $||u||=1$ so that
$T$ is a UGRF. Without loss of generality
we can assume that $n=1$ and hence $f(z) = z$. It is easy to see that
$R_t = [t,+\infty)$ and further
$$
  {\rm Tube}(R_t, r) = [t-r, +\infty) = R_{t-r}.
$$
This observation leads directly to the EC density of the Gaussian random field
\begin{equation}\label{rhoz}
\rho_d^{\rm G}(t)=\left(\frac{-1}{\sqrt{2\pi}}\frac{\partial}{\partial
t}\right)^d \Pp(T\geq t).
\end{equation}
We shall exploit this observation, that the tube is another rejection region but with a lower threshold, to derive the EC density for the $\bar\chi$ random field in the next section.

\subsection{The $\bar\chi$ random field}
\label{sec:bar}

Now let $R_t \subset \real^n$ be the rejection region for the $\bar\chi$ random field
at level $t$. This rejection region is the union of half planes all a distance $t$ from the origin. It is clear that a tube of radius $r$ about such a rejection region is simply another union of half planes all a distance $t-r$ from the origin (provided $r<t$). We thus arrive at precisely the same expression as for the Gaussian case:
${\rm Tube}(R_t, r) = R_{t-r}$.
In exactly the same way, this leads directly to the following representation for the
EC densities of a $\bar\chi$ random field:
\begin{equation}\label{rhochibar:pre}
\rho_d^{\rm \bar\chi}(t)=\left(\frac{-1}{\sqrt{2\pi}}\frac{\partial}{\partial
t}\right)^d \Pp(\bar\chi\geq t).
\end{equation}
We can now use the mixture representation (\ref{eq:Pchibar}) to show that the EC density of $\bar\chi$ is the same mixture of EC densities of the $\chi_j$ random field.
To see this, note that, by setting $U=\sphere{j}$ in (\ref{rhochibar:pre}), the EC density of $\chi_j$ is
\begin{equation}\label{eq:EC:chi}
\rho_d^\chi(t; j)=\left(\frac{-1}{\sqrt{2\pi}}\frac{\partial}{\partial
t}\right)^d \Pp(\chi_j\geq t).
\end{equation}
Combining this with (\ref{rhochibar:pre}) and (\ref{eq:Pchibar}) leads to the first expression of the following Theorem.

\begin{theorem}\label{thm:bc}
If $\cone(U)$ is convex then the EC density of the $\bar\chi$ random field is
$$
\rho_d^{\bar\chi}(t)=\sum_{j=1}^n p_j(U) \rho_{d}^{\chi}(t; j)=\sum_{j=0}^{n-1}
\lips_j(U) \rho^G_{d+j}(t)
$$
where $\rho_d^\chi(t; j)$ and $\rho_d^{\rm G}(t)$
are the EC densities of the the $\chi_j$ random field (\ref{eq:EC:chi}) and Gaussian random field (\ref{rhoz}), respectively.
\end{theorem}

The second part of the Theorem is proved as follows. Another way of evaluating $\Pp(\bar\chi\geq t)$ is to note that $u'Z$, as a function of $u$, is a UGRF and that $\bar\chi$ is its maximum over $U$. Hence we can
use the approximation \eqref{p} for Gaussian random fields, replacing $S$ by $U$. This is exact for $t > 0$ when $\cone(U)$ is convex. The reason is that the excursion set $\{u\in U : u'Z \geq t \}$ generates a cone that is the intersection of a convex circular cone (provided $t>0$) with convex $\cone(U)$, which is again convex. The EC of $\{u\in U : u'Z \geq t \}$ is either 0 or 1 if it is empty or not, that is, if $\bar\chi$ is less than or greater than $t$. Hence the expected EC is the P-value, so that \eqref{p} is exact and gives
\begin{equation}\label{rhochibar2}
\Pp(\bar\chi\geq t)=\sum_{j=0}^{n-1}
\lips_j(U) \rho^G_{j}(t).
\end{equation}
Combining this with (\ref{rhochibar:pre}) yields the second expression of Theorem \ref{thm:bc}. Note that the weights $p_j(U)$ can now be expressed in terms of
intrinsic volumes by equating (\ref{rhochibar2}) to (\ref{eq:Pchibar}) to give
$$
p_j(U)=\frac{1}{2^{j}\pi^\frac{j-1}{2}\Gamma(\frac{j+1}{2})}\sum_{m=0}^{\lfloor(n-j)/2\rfloor}
\frac{(-1)^m (d+2m)!}{(4\pi)^m m!}\lips_{j+2m-1}(U)
$$
(see Chapter 15 in \citet{RFG}).

\bigskip

\noindent{\bf Remark 1:} If $\cone(U)$ is not convex, the above argument used to derive \eqref{rhochibar2} fails, 
though \eqref{rhochibar:pre} still holds for the coefficients in the {\em exact} tube expansion, in the sense that ${\rm Tube}(R_t,r)=R_{t-r}$. However, if $\cone(U)$ is {\em locally convex} \eqref{rhochibar2} is exponentially accurate \citet{TTA:Validity} and therefore the right hand side of the result in Theorem 1 is the EC density up to an exponentially small error.

\bigskip

\noindent{\bf Remark 2:} The representation \eqref{eq:bar:chi:dec} represents
$\bar\chi(s)$ (reinstating dependence on $s$) as a mixture of $\chi_j(s)$
random fields with weights $p_j(U)$. It is therefore not surprising that
the EC density of the $\bar\chi(s)$ random field is a mixture of
the EC densities of $\chi_j(s)$ random fields with the
same weights. We give a sketch of a proof why this should be so for
the simplest cone: the positive orthant in $\real^k$
$$
\bar\chi(s)^2 =\sum_{j=1}^k 1_{\{Z_j(s) > 0\}} Z_j(s)^2.
$$
For this cone, a face is determined by a subset of $\{1, \dots, k\}$
which are the set of non-negative components of $\widehat{\mu}(s)$.
It is not hard to see that $\widehat{A}(s) = \{j: Z_j(s) < 0\}^c$ with the empty set
representing the vertex of the cone. We shall now make use of Morse theory, which shows that the EC of a set is determined by the critical points of a twice differentiable Morse function defined on the set \citep{Adler:1981}. The Morse theory expression for the EC density (\ref{rho}) is obtained by using the random field itself as the Morse function \citep{Worsley:1995:Boundary}.
The random field $\bar\chi(s)$ as a Morse function is actually differentiable (though not
twice differentiable) and it is not  hard to show that its critical
points are almost surely contained in the interior of the patches. This is
because the critical points on the boundary are points where a
particular
$\chi_j(s)$ random field has a critical point and one or more components
are 0 (see Figure \ref{fig0d0}).
For instance, critical points that appear on the segment of
boundary of the intersection of $\{s:Z_1(s)=0\}$ and the patch $\{s:\widehat{A}(s)=\emptyset\}$ are
points where $Z_1(s)$ has a critical point and $Z_1(s)=0$. The number of
such points is almost surely 0. Because there are no critical points
on the boundary of the patches, we can redefine $\bar\chi(s)$ near
these boundaries to get a Morse function with the same critical points
as $\bar\chi(s)$ and the standard Morse-theoretic
computation of the expected EC now shows that for each patch $J
\subset \{1, \dots, k\}$ we must find the number of critical points of
$\chi_J(s)^2 = \sum_{j \in J} Z_j^2(s)$ above the level $t$, counting
multiplicities. The expected EC above the level $t$, similar to (\ref{rho}) will therefore be
\begin{equation*}
  \begin{aligned}
    \sum_{J \subset \{1, \dots, N\}} \Ee \left(1_{\{\widehat{A}(s)=J\}}1_{\{\chi_J(s) > t\}}{\rm det}(- \ddot \chi_{J,d}(s))\ |\ \dot
\chi_{J,d}(s)=0\right){\mathbb P}(\dot \chi_{J,d}(s) = 0).
  \end{aligned}
\end{equation*}
Noting that the conditional distribution of $\ddot\chi_{J,d}(s)$ given $(Z(s),\dot Z(s))$ depends on $Z(s)$ only through $\|Z_J(s)\|$ implies that
$\ddot\chi_{J,d}(s)$ and $1_{\{\widehat{A}(s)=J\}}$ are conditionally independent given $(Z(s), \dot{Z}(s))$. In fact, this also implies that they are actually unconditionally independent.
This completes the
sketch of the proof: the sum over all subsets $J$ of size $j$
yields $p_j(U)$ times the EC densities of
$\chi^2_j$ random fields from (\ref{rho}). To go from the $\bar\chi(s)$ to the $T_{\rm
  IN}(s)$ or $T_{\rm LR}(s)$ random field is not complicated: simply replace
$\chi_J$ above by the appropriate $F$ random fields in the
decomposition \eqref{eq:LR:dec} or \eqref{eq:IN:dec}, though the
conditional independence argument is just slightly more complicated.
In the following sections, we prefer to use the Gaussian kinematic formula to
give a more direct and complete proof which
does not refer to Morse theory and counting critical points.

\subsection{The F- and T-statistic random fields}

\begin{figure}[t] \begin{center}
\mbox{\scalebox{0.5}{\includegraphics*{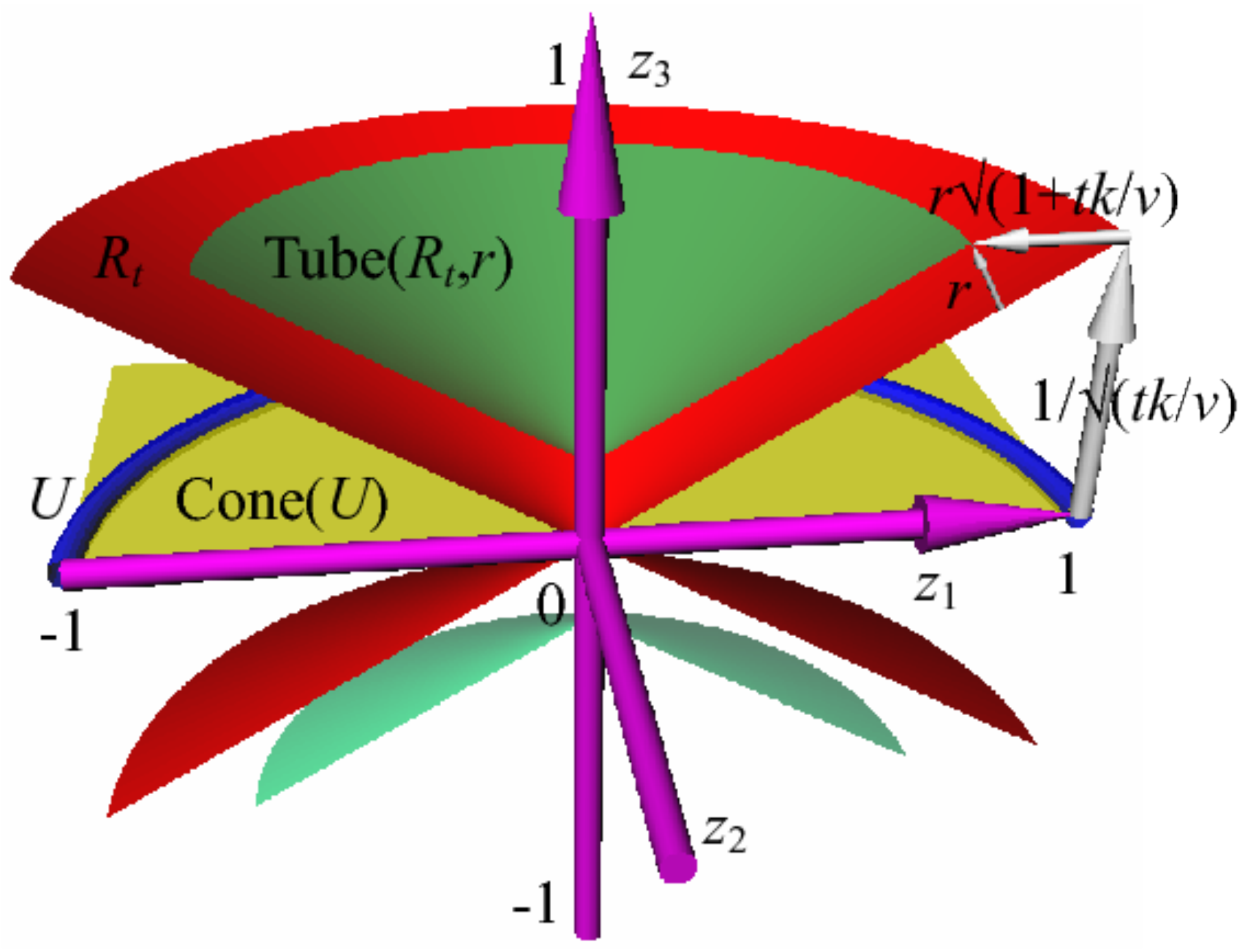}}}
\end{center}
\caption{Rejection region $R_t$ of the F statistic $F=(z_1^2+z_2^2)/2/z_3^2$ with $k=2$ and $\nu=1$. The purple axes are from
-1 to 1. The cone generator $U$ is blue, $\cone(U)$ is transparent yellow. The rejection region for a threshold
of $t=3/2$ is red; the tube about the rejection region (radius $r=0.15$) is
transparent green. Both rejection region and tube are cut at $z_2\geq 0$ and $|z_3|\leq 1/\sqrt{3}$. We expand the probability of this tube as a power series in $r$; its coefficients are the EC densities we seek.
}\label{fig:F}
\end{figure}

Our main results, stated in Theorem \ref{thm:IN} and Theorem \ref{thm:LR}, are based on a simple
refinement of Theorem \ref{thm:bc} in which we incorporate a $\chi^2$ field in the denominator. To see how it works, let us use the Gaussian kinematic formula to derive the EC density of the F-statistic field. Let
$R_{t} \subset \real^n$ be the rejection region of the F-statistic
random field $F$ with $k,\nu$ degrees of freedom. Without loss of generality, setting $z=(z_1,\dots,z_n)$, we can take
$$
f(z)=\frac{\sum_{i=1}^{k}z_i^2/k}{\sum_{i=k+1}^{n}z_i^2/\nu}.
$$
Then, a little elementary geometry (see Figure \ref{fig:F}) shows that
\begin{equation}\label{rhoF}
\Pp \left(Z\in {\rm Tube}(R_{t},r)\right)= \Pp
(\chi_k\geq T_r)+O(r^{n})
\end{equation}
where
$$
T_r=\chi_\nu\sqrt{\frac{tk}{\nu}}-r\sqrt{1+\frac{tk}{\nu}}.
$$
The remainder above reflects the fact that the tube ${\rm Tube}(R_{t}, r)$ is {\em almost} equal to the event $\{\chi_k\geq T_r \}$.
Near the origin, this fails but the probability content of where this fails is of order $O(r^{n})$. Further, the EC densities of $F$ are only defined for $d\leq D<n$ (as explained in Section \ref{sec:dim}).
Continuing with the main term in \eqref{rhoF}, and making use of (\ref{rhoz}),
\begin{equation}\label{eq:PchiTr}
  \begin{aligned}
\Pp
(\chi_k\geq T_r)
& = \Ee \left(\Pp \left(\chi_k\geq T_r \ \biggl| \chi_{\nu}\right) \right) \\
& = \Ee \left(\sum_{j=0}^{k-1} \lips_j(O^{k-1}) \; \rho^G_j( T_r) \right) \\
& = \sum_{d=0}^{\infty} \frac{(2\pi)^{d/2}r^d}{d!} \left(1+\frac{tk}{\nu} \right)^{d/2}\sum_{j=0}^{k-1} \lips_j(O^{k-1}) \; \Ee \left( \rho^G_{j+d}\left( \chi_\nu\sqrt{\frac{tk}{\nu}}\right) \right). \\
  \end{aligned}
\end{equation}
Hence, the EC densities for an F-statistic random field with $k,\nu$ degrees of freedom are given by
\begin{equation}
  \label{eq:rhoF:full}
\rho^{\rm F}_d(t; k, \nu) =  \left(1+\frac{tk}{\nu} \right)^{d/2}\sum_{j=0}^{k-1} \lips_j(O^{k-1}) \; \Ee \left( \rho^G_{j+d}\left( \chi_\nu\sqrt{\frac{tk}{\nu}}\right) \right).
\end{equation}

For the T-statistic random field $T_1$, a similar argument to that leading to \eqref{rhoF} shows that we must expand the following probability in a power series:
$$
\Pp
\left(Z_1\geq \chi_\nu\sqrt{\frac{t^2}{\nu}}-r\sqrt{1+\frac{t^2}{\nu}}\right)
$$
where $Z_1 \sim N(0,1)$ is independent of $\chi_{\nu}$. In the above expression, $t^2$ appears instead of $t$ because $T_1^2$ is an $F_{1,\nu}$ random field and $Z_1$ appears rather than $\chi_1=|Z_1|$ on the left hand of the inequality side because $T_1$ is one-sided.
Similar calculations to those above for the F-statistic yield the following expression for the EC densities of the T-statstic random field
$$
  \begin{aligned}
\rho^{\rm T}_d(t;\nu)& = \left(1+\frac{t^2}{\nu} \right)^{d/2} \Ee \left( \rho^G_{d}\left( \chi_\nu\sqrt{\frac{t^2}{\nu}}\right) \right)\\
 &= \sum_{l=0}^{\lfloor \frac{d-1}{2} \rfloor} \frac{(-1)^l(d-1)! \Gamma\left(\frac{d-1-2l+\nu}{2}\right)}{\pi^{(d+1)/2}2^{2l+1}(d-1-2l)! l! \Gamma\left(\frac{\nu}{2}\right)}
 \left(\frac{t^2}{\nu}\right)^{(d-1-2l)/2} \left(1 + \frac{t^2}{\nu} \right)^{-(\nu-1-2l)/2}
\end{aligned}
$$
for $d>0$ and $\mathbb{P}(T_1>t)$ for $d=0$. This is simpler than the expression in \citet{Worsley:1994:Chi:t:F}; it is a single sum, whereas the the expression in \citet{Worsley:1994:Chi:t:F} is a double sum.

A simple rearrangement of \eqref{eq:rhoF:full} yields the following equivalent representation of  the EC densities of the F-statistic random field in terms of the EC densities of the T-statstic random field:
$$
\rho^{\rm F}_d(t; k, \nu) =  \left(1+\frac{tk}{\nu} \right)^{-d/2}\sum_{j=0}^{k-1} \lips_j(O^{k-1}) \; \rho^{\rm T}_{d+j}(\sqrt{tk}; \nu).
$$

\subsection{The independently normalized cone random field $T_{\rm IN}$}\label{sec:GKF:IN}

\begin{figure}[t] \begin{center}
\mbox{\scalebox{0.5}{\includegraphics*{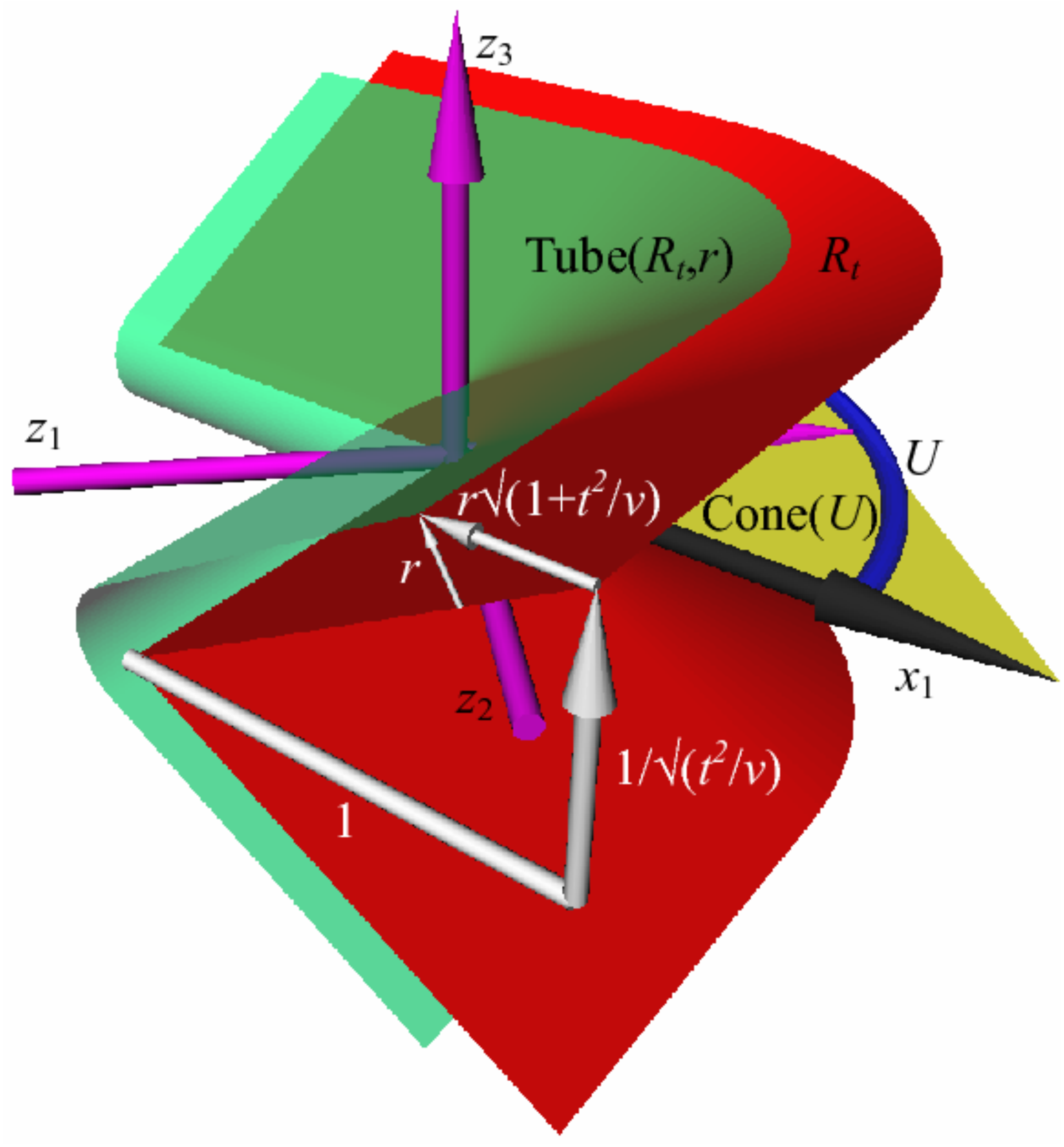}}}
\end{center}
\caption{Rejection region $R_t$ of the independently normalized test statistic $T_{\rm IN}$ for the same cone as in Figure \ref{fig0d0} and the same $z$ as in Figure \ref{fig:F}. The cone edges $x_1$ and $x_2$ are black. The threshold
is $t=\sqrt{3}$ and both the rejection region and tube are cut at $z_1\pm z_2\geq -\sqrt{2}$ and $|z_3|\leq 1/\sqrt{3}$.
}\label{fig:coneIN}
\end{figure}

It is slightly easier to work with $T_{\rm IN}$, since it more closely resembles $F$, so we tackle this ahead of $T_{\rm LR}$. It should now be clear how to proceed:
find the rejection region as a function of the $n$ UGRF's; put a tube around
with radius $r$; work out the probability content; differentiate $d$ times to
get the EC density. This sounds formidable, but it is in fact virtually identical
to the case of the F-statistic presented above.
For readers with good
geometric intuition, Figure \ref{fig:coneIN} might help: it shows the simple case of the rejection region $R_t=\{ Z: T_{\rm IN}\geq t\}$ where $k=2$
and $\nu=1$, and $U$ is a quarter circle, as in Figure \ref{fig0d0}.

\begin{theorem}\label{thm:IN}
If $\cone(U)$ is convex then the EC density of the independently normalized cone random field $T_{\rm IN}$ is
$$
\rho_d^{\rm IN}(t)=
\sum_{j=1}^{k} p_j(U) \rho_d^{\rm F}\left(\frac{t^2}{j};j,\nu\right)=
\sum_{j=0}^{k-1} \lips_j(U) \rho_{d+j}^{\rm T}(t;\nu)
\left(1+\frac{t^2}{\nu}\right)^{-j/2}.
$$
The EC densities are valid
for $d < \nu + \max(l(U),1)$, where $l(U)$ is the dimension of the largest linear subspace in $\cone(U)$.

\end{theorem}

\noindent{\bf Remark:}  The representation \eqref{eq:IN:dec} represents
$T_{\rm IN}$ as a patchwork mixture of $\sqrt{j \cdot F_{j,\nu}}$
random fields with weights $p_j(U)$.  See Remark 2 after Theorem \ref{thm:bc} for why Theorem \ref{thm:IN} should not be surprising.
For the case of non-convex $\cone(U)$, see Remark 1 after Theorem \ref{thm:bc}.

\bigskip

\noindent{\bf Proof:}
The same geometric argument that led to \eqref{rhoF} leads to the following
approximate equality
$$
\left\{ Z \in {\rm Tube}(R_t,r) \right \}\simeq\left\{ \bar\chi \geq
T_r^*\right\} \\
$$
where
$$
T_r^* = \chi_\nu \sqrt{\frac{t^2}{\nu}} - r \sqrt{1 + \frac{t^2}{\nu}}.$$
In fact, $\{Z \in {\rm Tube}(R_t,r)\}$ is contained within $\{\bar\chi \geq T_r^*\}$ with the difference
coming from points where $T_r^*$ and $\bar\chi$ are both near 0.  If $l(U) > 1$, the probability of this difference, as a function of the tube radius $r$, is of order $O(r^{l(U)+\nu})$. If $l(U)=0$, then similar arguments to those in Section \ref{sec:dim} show that we need only worry about 0/0 when $\bar\chi>0$ but is close to 0, that is, when its $\chi_1$ components are near 0 and $\chi_{\nu}$ is also near 0. The probability of this is of order $O(r^{\nu+1})$.
Since we must have  $d<\nu+\max(l(U),1)$ anyway to avoid $0/0$, we can ignore this difference in either case, thus for our
purposes we need only expand ${\mathbb P}( \bar\chi \geq  T_r^*)$ as a power series in $r$. This computation is essentially identical
to the case of the F-statistic where $O^{k-1}$ is replaced with a general $U$.
Following the calculations preceding \eqref{eq:rhoF:full}:
$$
\begin{aligned}
{\mathbb P}( \bar\chi \geq
T_r^*)&={\mathbb E}\left( \sum_{j=0}^{k-1}\lips_j(U) \rho_j^{\rm
G}(T_r^*)\right) \\
&=  \sum_{d=0}^{\infty} \frac{(2\pi)^{d/2}r^d}{d!} \left(1+\frac{t^2}{\nu} \right)^{d/2}\sum_{j=0}^{k-1} \lips_j(U) \; \Ee \left( \rho^G_{j+d}\left( \chi_\nu\sqrt{\frac{t^2}{\nu}}\right) \right) \\
&=  \sum_{d=0}^{\infty} \frac{(2\pi)^{d/2}r^d}{d!} \sum_{j=0}^{k-1}  \lips_j(U) \; \rho^{\rm T}_{j+d}(t; \nu)\left(1+\frac{t^2}{\nu} \right)^{-j/2}. \\
\end{aligned}
$$
To derive the EC densities in terms of $F$ EC densities, simply use (\ref{eq:Pchibar}), (\ref{eq:PchiTr}) and (\ref{eq:rhoF:full}):
$$
\begin{aligned}
{\mathbb P}( \bar\chi \geq
T_r^*)&= \sum_{j=\max(l(U),1)}^{k} p_j(U) \; \Pp\left(\chi_j \geq T_r^* \right)  \\
&= \sum_{d=0}^{\infty} \frac{(2\pi)^{d/2}r^d}{d!} \sum_{j=\max(l(U),1)}^{k} p_j(U) \; \rho^{\rm F}_d \left(\frac{t^2}{j}; j, \nu\right)
\end{aligned}
$$
\qed

\subsection{The likelihood ratio cone random field $T_{\rm LR}$} \label{sec:GKF:LR}

\begin{figure}[t] \begin{center}
\mbox{\scalebox{0.5}{\includegraphics*{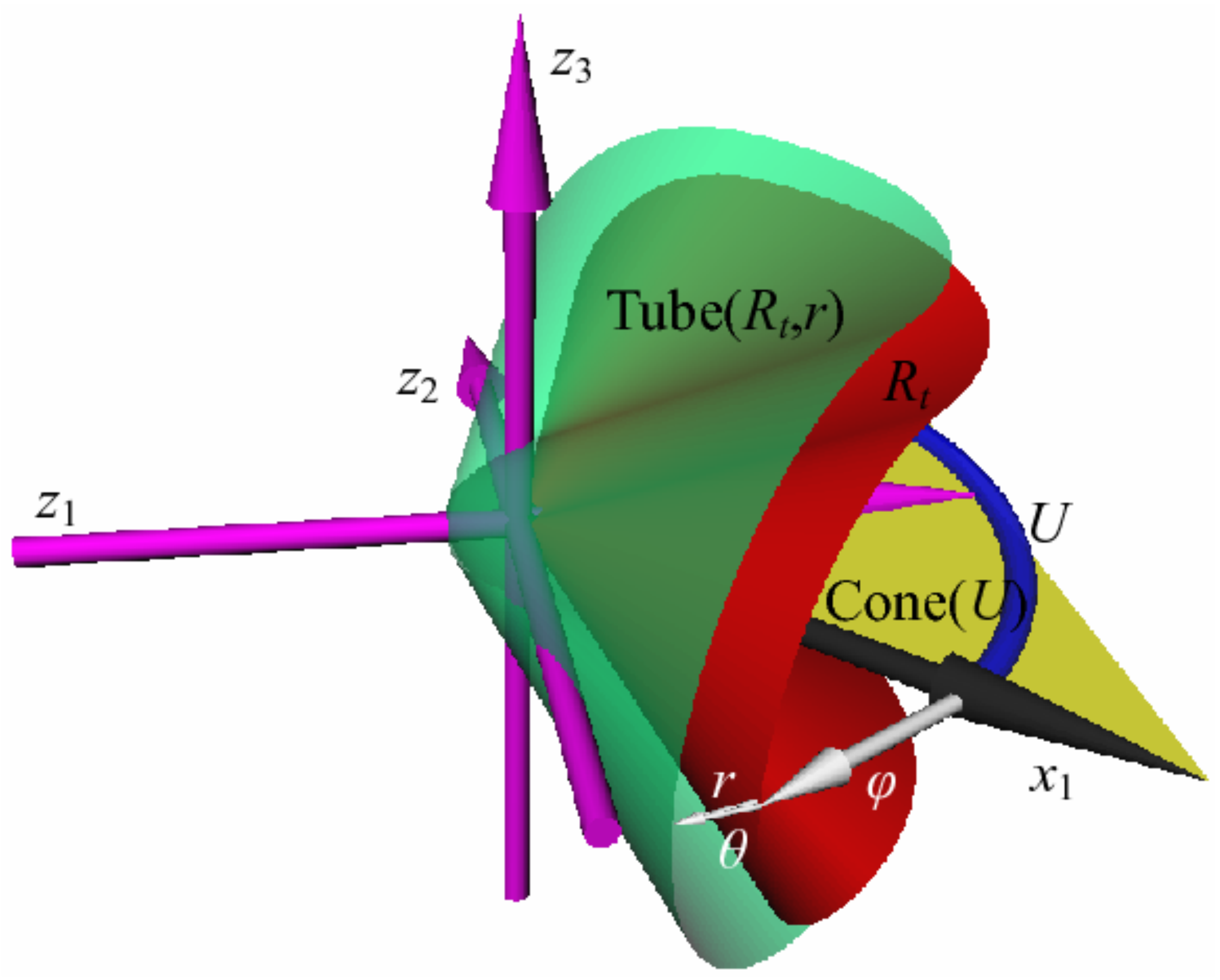}}}
\end{center}
\caption{As for Figure \ref{fig:coneIN}, but for the likelihood ratio test statistic $T_{\rm LR}$ at a threshold $t=3$, cut at $||z||\leq 1$; $\phi=\arccos(t/\sqrt{n+t^2})=\pi/6$.
}\label{fig:coneLR}
\end{figure}

Figure \ref{fig:coneLR} illustrates the rejection region $R_t$ of $T_{\rm LR}$.

\begin{theorem}\label{thm:LR}
If $\cone(U)$ is convex then the EC density of the likelihood ratio cone random field $T_{\rm LR}$ is
$$
\rho_d^{\rm LR}(t)=\sum_{j=1}^{n} p_j(U) \rho_d^{\rm
F}\left(\frac{t^2}{j}\frac{n-j}{n};j,n-j\right)
$$
The EC densities are valid for $d < n$.
\end{theorem}

\noindent{\bf Remark:} As for $T_{\rm IN}$, the representation \eqref{eq:LR:dec} represents
$T_{\rm LR}$ as a patchwork mixture of $\sqrt{jn/(n-j) \cdot F_{j,n-j}}$
random fields with weights $p_j(U)$. See Remark 2 after Theorem \ref{thm:bc} for why Theorem \ref{thm:LR} should not be surprising. For the case of non-convex $\cone(U)$, see Remark 1 after Theorem \ref{thm:bc}.

\bigskip

\noindent{\bf Proof:}
It is easier to transform to the equivalent correlation coefficient
$$
C=\frac{T_{\rm LR}}{\sqrt{n+T_{\rm LR}^2}}=\frac{\bar\chi}{||Z||}=\max_{u\in
U}\frac{u'Z}{||Z||}.
$$
Then the rejection region $C\geq c$ is simply a cone centered at the origin that
intersects the unit sphere in a tube of geodesic radius $\phi=\arccos
c=\arccos(t/\sqrt{n+t^2})$ about $U$:
$$
R_t=\left\{z: \arccos\left( \max_{u\in U}\frac{u'z}{||z||}\right)\leq
\phi\right\}.
$$
When $\cone(U)$ is convex there is an exact expression for the probability content of a tube
about a subset of the sphere, similar to (\ref{eq:Pchibar})
\citep{Lin:Lindsay:1997,Takemura:Kuriki:Chi:Bar}:
$$
{\mathbb P}\left(\frac{\bar\chi}{||Z||}\geq c\right)={\mathbb P}(Z\in R_t)=
\sum_{j=1}^{n} p_j(U) {\mathbb P}\left( \arccos(\sqrt{B_{j}})\leq \phi \right)
$$
where $B_j$ is a Beta random variable with parameters $j/2,(n-j)/2$ (with
$B_{n}=1$ with probability one). The restriction of $\cone(U)$ to a convex set is not necessary, as it was for $\bar\chi$ - the only requirement is
that $t$ must be sufficiently large (i.e. $\phi$ must be sufficiently small) so
that the tube does not self-intersect. This phenomenon is similar to what occurs
when establishing the  accuracy of \eqref{p} for non-convex
regions $\cone(U)$.
If $\cone(U)$ is convex then $t\geq 0$
suffices.

The next step is to put a tube about the rejection region $R_t$. Provided $r$ is
sufficiently small, a (Euclidean) tube of radius $r$ about $R_t$ intersects the
sphere of radius $||z||$ in a spherical tube of geodesic radius
$\theta=\arcsin(r/||z||)$ about $R_t$. For fixed $||z||$ sufficiently large,  $R_t$ is already a spherical tube about $||z||U$, so the (Euclidean) tube about $R_t$ is a spherical tube about $||z||U$ of geodesic radius $\phi+\theta$:
$$
{\rm Tube}(R_t,r)=\left\{z: \arccos\left( \max_{u\in
U}\frac{u'z}{||z||}\right)\leq \phi +\theta \right\}.
$$
The part of the tube near the origin with small $||z||$ may
contain a ``wedge'' of the ball of radius $r$ (see Figure \ref{fig:coneIN}(a)) that is the
only part of the whole tube that contributes to the coefficient of $r^n$.
 As pointed out in
Section \ref{sec:dim}, $T_{\rm LR}$ is only defined for $d\leq D<n$
so we can ignore this.
 It therefore follows
that it is sufficient for us to work with
\begin{equation}\label{eq:B}
{\mathbb P}\left( Z \in {\rm Tube}(R_t,r) \right)=\sum_{j=1}^{n} p_j(U) {\mathbb
P}\left( \arccos(\sqrt{B_{j}}) \leq \phi +\Theta \right)+O(r^n),
\end{equation}
where $\Theta=\arcsin(r/||Z||)$ is independent of $B_{j}$. The inequality in
(\ref{eq:B}) is
$$
\arccos(\sqrt{B_{j}})-\phi \leq \Theta \Longleftrightarrow \sqrt{1-B_{j}} c
-\sqrt{B_{j}} \sqrt{1-c^2} \leq \frac{r}{||Z||},
$$
so that
$$
{\mathbb P}\left( \arccos(\sqrt{B_{j}}) \leq \phi +\Theta\right)=
 {\mathbb P}\left(\chi_{j}  \geq \chi_{n-j}\sqrt{ \frac{t^2}{n}} -
 r\sqrt{1+\frac{t^2}{n}}\right),
$$
where $\chi_{j}$ and $\chi_{n-j}$ are the square roots of independent $\chi^2$
random variables with degrees of freedom indicated by their subscripts. Putting
everything together, the EC density that we seek is the coefficient of
$r^d(2\pi)^{d/2}/d!$ in
$$
{\mathbb P}\left( Z \in {\rm Tube}(R_t,r) \right)=\sum_{j=1}^{n} p_j(U) {\mathbb
P}\left(\chi_{j}  \geq \chi_{n-j} \sqrt{\frac{t^2}{n}} -
 r\sqrt{1+\frac{t^2}{n}}\right)+O(r^n).
$$
Since this expression is linear in the tube probabilities, we can differentiate
immediately to arrive at the result we are looking for.
\qed

\section{Application}\label{sec:app}

\begin{figure}[t]
\begin{center}
\mbox{\scalebox{0.75}{\includegraphics*{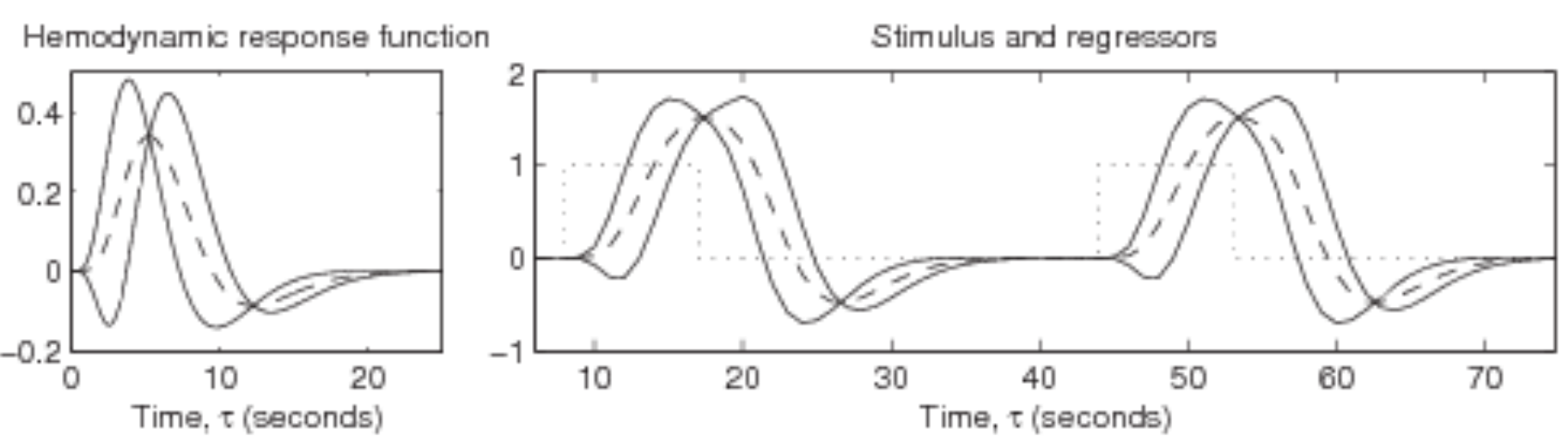}}}
\end{center}
\caption{The hemodynamic response function $h_0$ (left, dashed line) and the two extremes $h_0\pm 2\dot h_0$ (left, solid lines) convolved with the on-off painful heat stimulus $g$ (right, dotted line) to give the ``middle'' of the cone $u$ (right, dashed line) and the two cone edges, the regressors $x_{1,2}=(h_0\pm 2\dot h_0)\star g$ (right, solid lines). The on-off stimulus is repeated ten times, from 0 to 360 seconds. }\label{fighrf}
\end{figure}

\citet{Friman:2003} and \citet{Calhoun:2004} proposed the cone and one-sided F-statistics for the detection of functional magnetic
resonance (fMRI) activation in the presence of unknown delay in the hemodynamic response. We illustrate our methods with a re-analysis of the fMRI data from study an pain perception that was used by \citet{Worsley:Taylor:2005:latency}. The data, fully described in \citet{Worsley:fmristat}, consists of a time series of 3D fMRI
images $Z(s,\tau)$ at point $s\in \real^3$ in the brain at time $\tau$. The subject received an alternating 9
second painful then neutral heat stimulus to the right calf, interspersed with
9 seconds of rest, repeated 10 times. The mean of the fMRI data is modeled as the
indicator for each stimulus ($g(\tau)=1$ if on, 0 if not) convolved with a
known hemodynamic response function (hrf) $h_0(\tau)$ that delays and disperses the
stimulus by about 5.5 seconds (see Figure \ref{fighrf}). Taking $g(\tau)$ as just the painful heat stimulus, we add this to a linear model for the fMRI data:
$$
Z(s,\tau)=(h_0 \star g)(\tau) \beta(s) + \sigma(s)\epsilon(s,\tau),
$$
where $\epsilon(s,\tau)\sim{\rm N}(0,1)$. Our main interest is to detect
regions of the brain that are `activated' by the hot stimulus, that is, points $s$ where $\beta(s)>0$.

There is often some doubt about the 5.5 second delay of the hrf, so to allow
for unknown delay, we shift $h_0(\tau)$ by an amount $\delta(s)$ and add $\delta(s)$ as
a parameter to the hrf. To keep the linear model, we then approximate the
shifted hrf by a Taylor series expansion in $\delta(s)$ \citep{Friston:1998}:
$$
h(\tau; \delta(s))=h_0(\tau-\delta(s))\approx h_0(\tau)-\delta(s)\dot h_0(\tau).
$$
The convolution of $h(\tau; \delta(s))$ with the stimulus $g(\tau)$ is then roughly
equivalent to adding the convolution of $-\dot h_0(\tau)$ with the stimulus as an extra regressor to give the linear model:
$$
Z(s,\tau)=(h_0 \star g)(\tau)\beta(s)-(\dot h_0 \star g)(\tau)\beta(s)\delta(s)+\sigma(s)\epsilon(s,\tau).
$$
However the key ingredient in the model is that there is some
structure to the coefficients dictated by the physical nature of the
regressors. It is strongly suspected that $\beta(s)>0$ and the shift is restricted to a range of known plausible values $\delta(s)\in[\Delta_1,\Delta_2]$.
In our example, we take  $[\Delta_1,\Delta_2]=[-2, 2]$ seconds. It is easy to see that the restrictions specify a non-negative-coefficient regression model
$$
Z(s,\tau)=x_1(\tau)\beta_1(s)+x_2(\tau)\beta_2(s)+\sigma(s)\epsilon_i(s,\tau),\ \ \ \ \ \beta_1(s)\geq 0, \beta_2(s)\geq 0,
$$
with regressors $x_j= (h -\Delta_j\dot h )\star g$, $j=1,2$,
illustrated in Figure \ref{fighrf}. The model is sampled at $n$ equal intervals over time and suppose for simplicity that the resulting observations are independent. Replacing dependence on $\tau$ by vectors in $\real^n$, the linear model is the same as (\ref{eq:linmod}) with $m=2$:
\begin{equation}\label{mod3}
Z(s)=x_{1}\beta_1(s)+x_{2}\beta_2(s)+\sigma(s)\epsilon(s),\ \ \ \ \ \beta_1(s)\geq 0, \beta_2(s)\geq 0,
\end{equation}
where $\epsilon(s)$ is a vector of $n$ iid stationary Gaussian random fields. This model (\ref{mod3}) is of course a 2D ($k=2$) cone alternative with cone angle
\begin{equation}\label{theta}
\alpha=\arccos\left( x_1'x_2/(||x_1||\cdot||x_2||)\right).
\end{equation}
The cone intrinsic volumes are $\lips_{0,1}(U)=1,\alpha$, and the $\bar\chi$ weights are $p_{1,2}(U)=1/2,\alpha/(2\pi)$. The ``middle'' of the cone is $u=(x_1+x_2)/2$, appropriately normalized, which of course corresponds to the unshifted model with $\delta=0$.

In practice our observations were temporally correlated and we added regressors
to allow for the neutral heat stimulus and a cubic polynomial in the scan time to allow
for drift, leaving $n=112$ effectively independent observations sampled every 3 seconds. The resulting $\alpha$, found by whitening the regressors and removing the effect of the added nuisance regressors before calculating (\ref{theta}), now depends on $s$ since the temporal correlation depends on $s$. However $\alpha$ was remarkably constant across the brain, averaging at $\alpha=1.06\pm 0.03$ radians or $60.9\pm 1.7^\circ$, so we take it as fixed at its mean value.

The search region $S$ is the entire brain. The error random fields $\epsilon_i(s)$
are not isotropic, so we must use Lipschitz-Killing curvatures of $S$ instead
of intrinsic volumes. The highest order term with $d=D$ makes the largest
contribution to the P-value approximation (\ref{p}), and fortunately there is a
very simple unbiased estimator for $\lips_D(S)$
\citep{Worsley:Nonisotropic:1999,TW:noniso}. At a particular voxel, let $E$ be
the $n\times 1$ vector of least-squares residuals from
(\ref{mod3}), and let $N=E/||E||$. Let $Q$ be the $n\times D$ matrix of
their spatial nearest neighbor differences, that is, column $d$ of $Q$ is
$N(s_2)-N(s_1)$ where $s_1,s_2$ are neighbors on lattice axis $d$. Then the
estimator of $\lips_D(S)$ is
$$
\widehat\lips_D(S)=\sum {\rm det}(Q'Q)^{1/2},
$$
where summation is taken over all voxels inside $S$
\citep{Worsley:Nonisotropic:1999,TW:noniso}. The result is
$\widehat\lips_3(S)=8086$, which is of course unitless. The lower order
Lipschitz-Killing curvatures are much more difficult to estimate, but they can
be very accurately approximated by those of a ball with the same volume, that is with radius $r=12.5$, to give $\widehat\lips_{0,1,2}(S)=1, 4\pi r, 2\pi r^2$.

\begin{figure}[t]
\begin{center}
(a)\mbox{\scalebox{0.28}{\includegraphics*{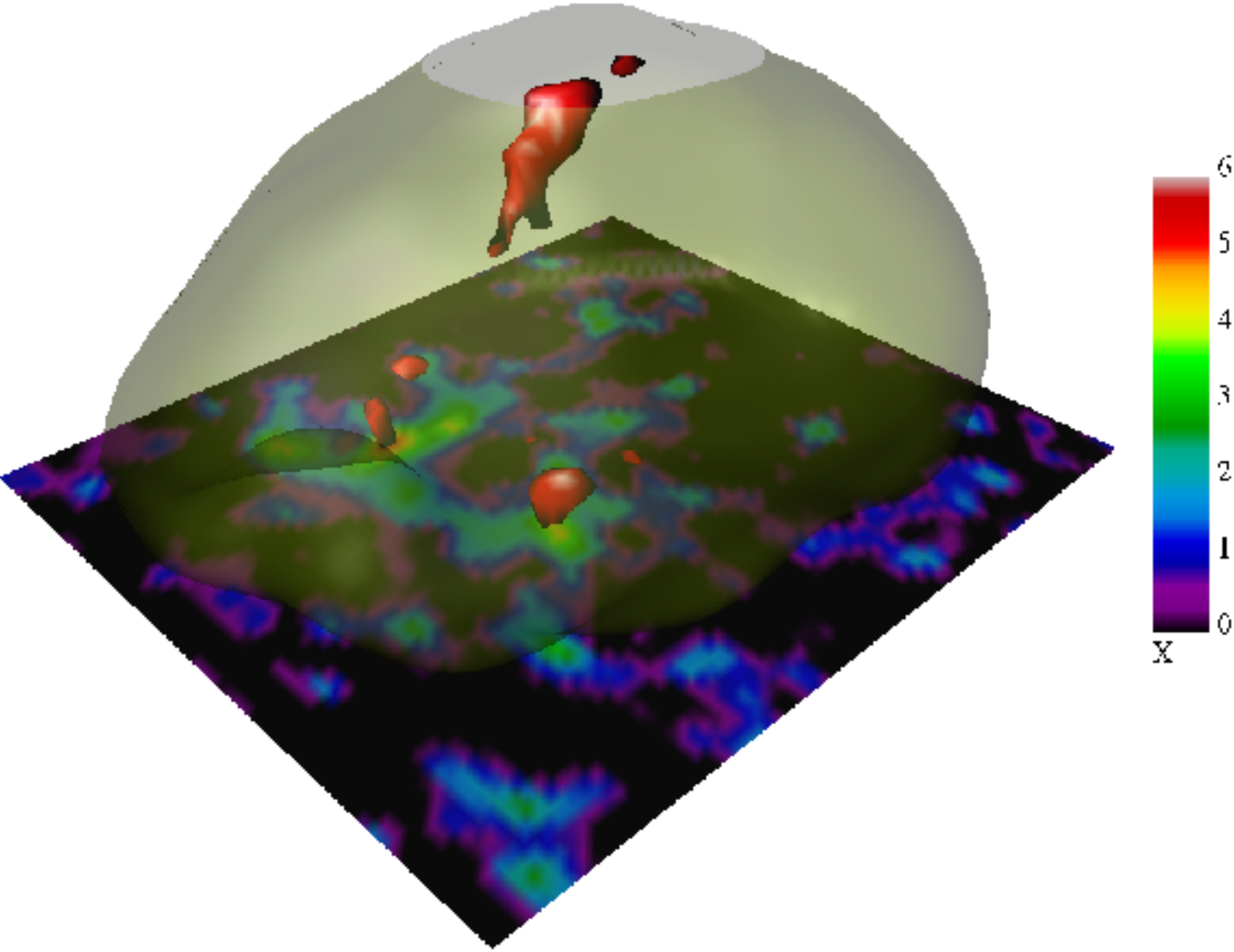}}}
(b)\mbox{\scalebox{0.28}{\includegraphics*{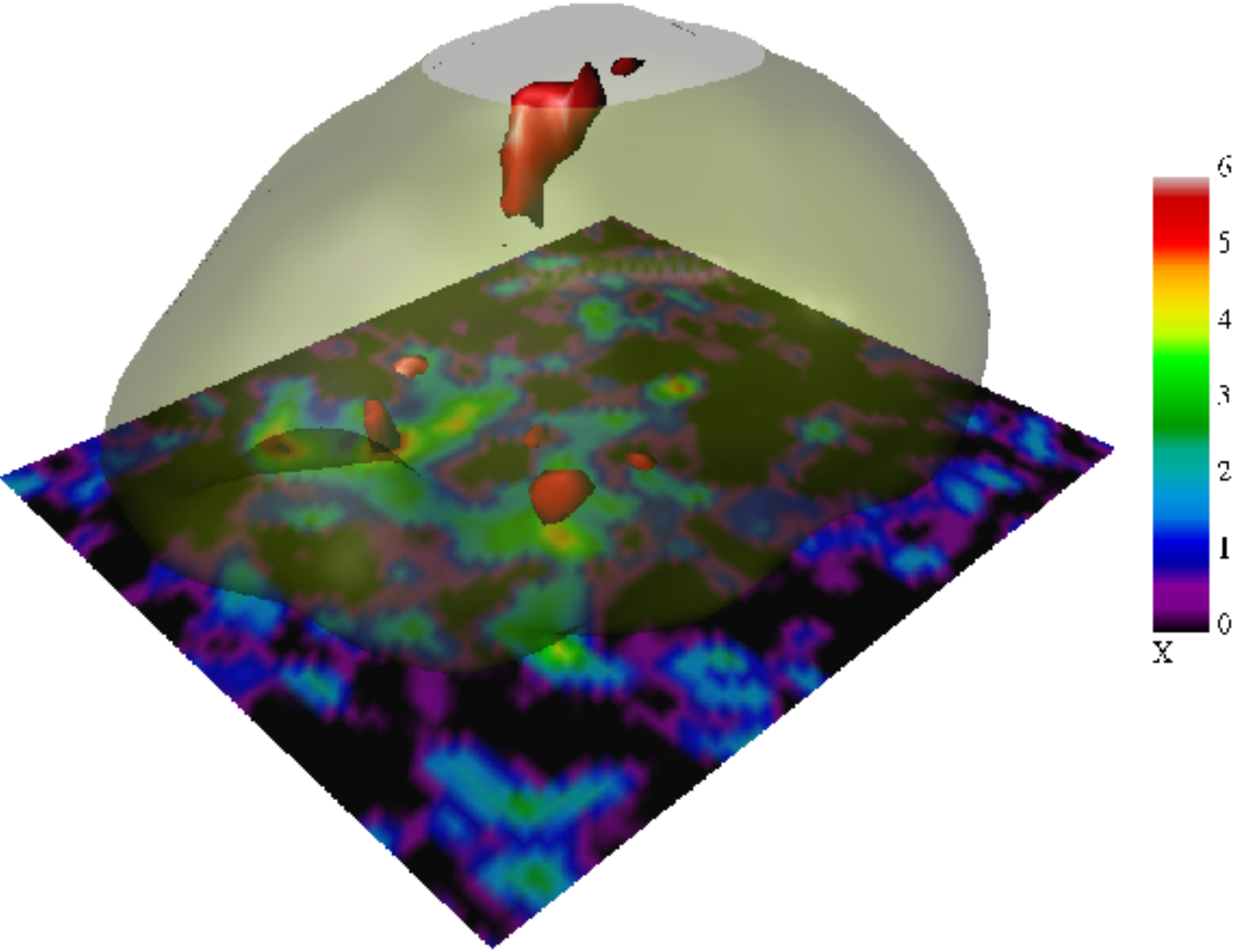}}}\\
(c)\mbox{\scalebox{0.28}{\includegraphics*{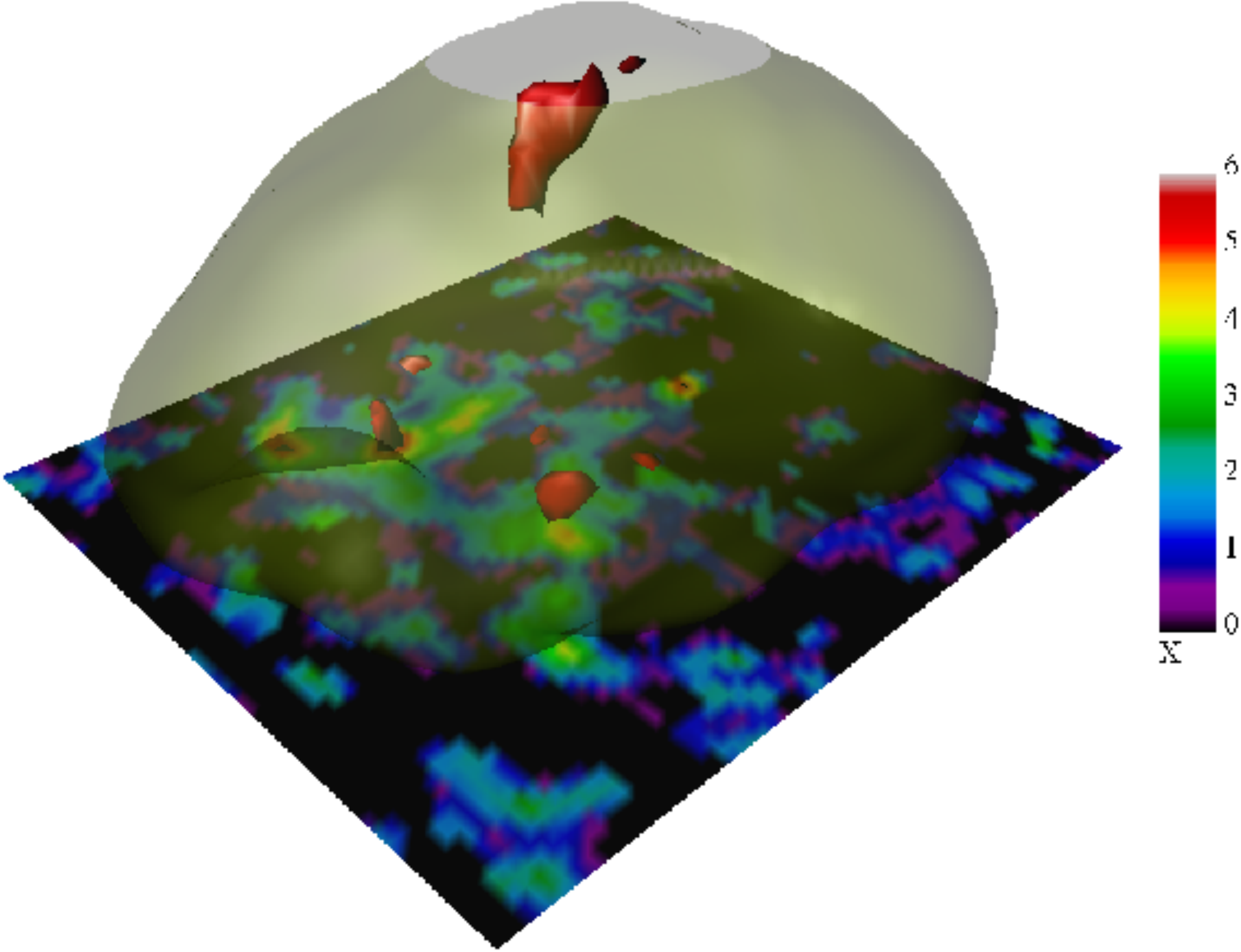}}}
(d)\mbox{\scalebox{0.28}{\includegraphics*{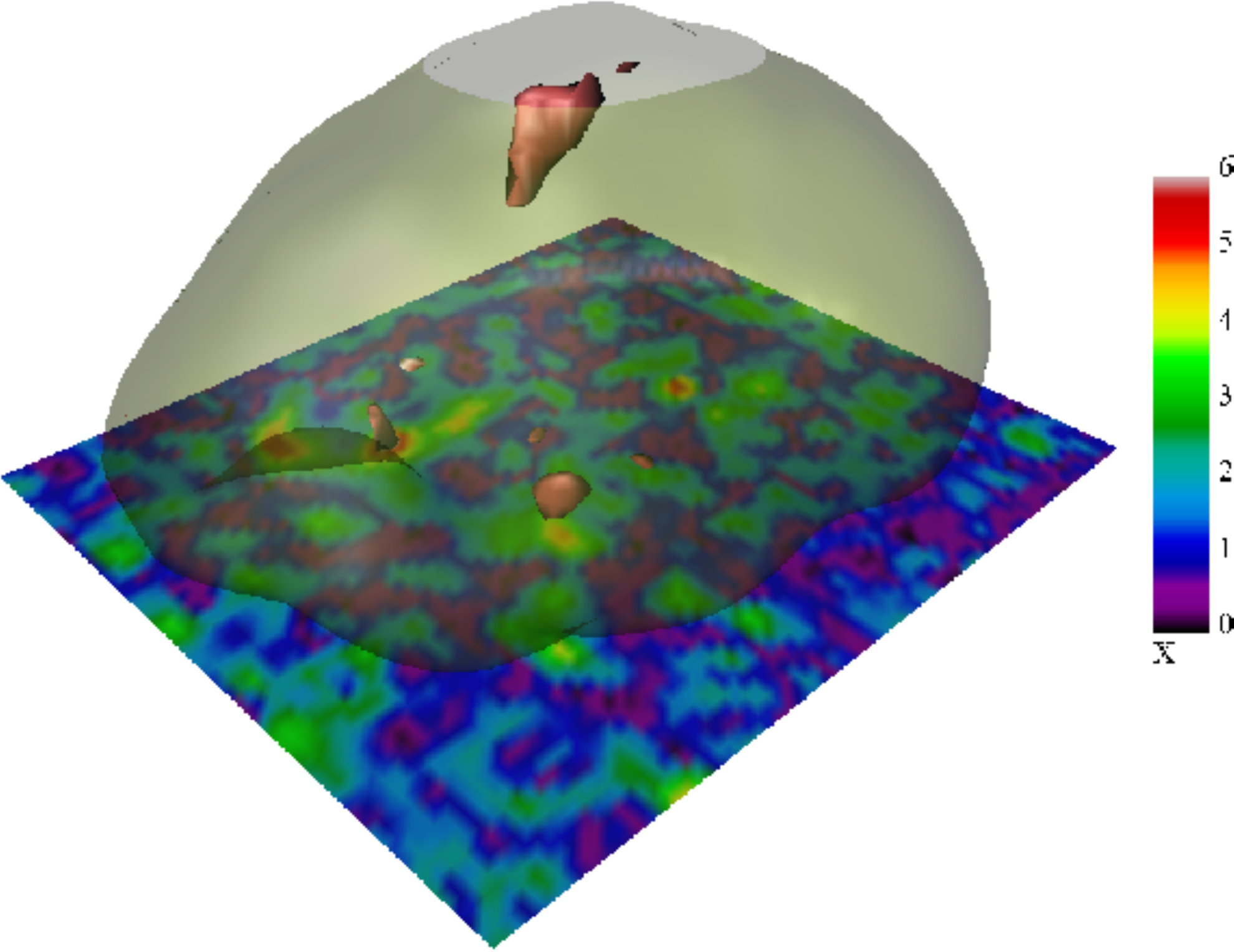}}}
\end{center}
\caption{Detecting activation in fMRI data. Each image shows the
search region (the brain, left front facing viewer) and a slice of the test
statistic (color coded) thresholded at $P=0.05$ (red-pink blobs - see Table
\ref{tab1}). The test statistics, in order of increasing threshold, are (a) the T-statistic $T_1$; (b) the
cone statistic $T_{\rm IN}$ (indistinguishable from $T_{\rm LR}$ in this case); (c) the square root of twice the
one-sided F-statistic $\sqrt{2F_+}$; (d) the square root of twice the
F-statistic $\sqrt{2F}$.}\label{figtstat}
\end{figure}

We are now ready to use (\ref{p}) to get approximate P-values for the maximum
of our test statistic random fields. Since the
degrees of freedom $\nu=110$ is so large, the two cone statistics were almost
identical, so we only show results for the independently normalized cone statistic. The
$P=0.05$ thresholds are shown in Table \ref{tab1}. Note that the values of the
statistics are increasing since the cone is getting larger, but of course the
$P=0.05$ thresholds are increasing as well to compensate for this. The net
result is that the volume of detected activation due to the painful heat stimulus
remains roughly the same. Interestingly, it is the cone statistic
with delays in the range $[-2,2]$ seconds that detects the most activation.
This activation is shown in Figure \ref{figtstat} (left primary somatosensory
area and left and right thalamus).

The last question is which test is the most powerful.
\citet{Worsley:Taylor:2005:latency} gives a power comparison of the four
tests that shows that if the true delay is in the range $[-1,1]$ seconds then
the usual T-statistic $T_1$ is the most powerful, but outside this range, the
cone statistic is the most powerful.

\begin{table}\begin{center}
\begin{tabular}{lcc}
Test statistic & $P=0.05$ threshold & Detected volume (cc) \\
\hline
(a) T-statistic, $T_1$ & 5.15 & 4.0 \\
(b) Cone statistic, $T_{\rm LR}\approx T_{\rm IN}$ & 5.44 & 4.3\\
(c) One-sided F-statistic, $\sqrt{2F_+}$ & 5.63 & 3.8\\
(d) F-statistic, $\sqrt{2F}$ & 5.80 & 2.9\\
\end{tabular}\end{center}
\caption{Test statistics, $P=0.05$ thresholds, and volume of detected
activation for the application in Figure \ref{figtstat}, in order of increasing
threshold. The cone statistic detects the most
activation.}\label{tab1}\end{table}

\appendix

\section{Intrinsic volume}\label{app:iv}

The $d$-dimensional intrinsic volume of a set $S$ is a generalization of its
volume to lower dimensional measures. The $D$-dimensional intrinsic volume of
$S\subset \real^D$ is its usual volume or Lebesgue measure, the
$(D-1)$-dimensional intrinsic volume of $S$ is half its surface area, and the
$0$-dimensional intrinsic volume is the Euler characteristic of $S$. The
simplest definition is {\em implicit}, identifying the intrinsic volumes as
coefficients in a certain polynomial. This definition comes from the
Steiner-Weyl volume of tubes formula which states that if $S$ has no concave `corners', then for $r$ small enough
\begin{equation}\label{eq:voltub}
|{\rm Tube}(S,r)| = \sum_{d=0}^D \omega_{D-d} r^{D-d} \lips_d(S)
\end{equation}
where $|\cdot|$ denotes Lebesgue measure and $\omega_d= \pi^{d/2} /
\Gamma(d/2+1)$ is the Lebesgue measure of the unit ball in $\real^d$.

If $S$ is bounded  by a smooth hypersurface, so that there is a unique
normal vector at each point on the boundary, then a more direct definition is as follows.
Let $C(s)$ be the $(D-1) \times (D-1)$ inside curvature matrix at $s \in
\partial S$, the boundary of $S$. To compute the intrinsic volumes, we need the
{\em det-traces} of a square matrix: for a $d\times d$ symmetric matrix $A$,
let ${\rm detr}_j(A)$ denote the sum of the determinants of all $j\times j$
principal minors of $A$, so that ${\rm detr}_d(A)=\det(A)$, ${\rm
detr}_1(A)={\rm tr}(A)$, and we define ${\rm detr}_0(A)=1$. Let $a_d=2
\pi^{d/2}/\Gamma(d/2)$ be the $(d-1)$-dimensional Hausdorff (surface) measure
of the unit $(d-1)$-sphere in $\real^d$. For $d=0,\dots,D-1$ the $d$-dimensional
intrinsic volume of $S$ is
$$
\lips_d(S)={{1}\over {a_{D-d}}}\int_{\partial S} {\rm detr}_{D-1-d}\{C(s)\}ds,
$$
and $\lips_D(S)=|S|$, the Lebesgue measure of $S$. Note that
$\lips_0(S)=\varphi(S)$ by the Gauss-Bonnet Theorem, and $\lips_{D-1}(S)$ is
half the surface area of $S$.

For the unit $(k-1)$-sphere, $C=\pm I_{(k-1)\times (k-1)}$ on the outside/inside of $\sphere{k}$, so that
\begin{equation}\label{intvolsphere}
\lips_d(\sphere{k})= 2 {k-1 \choose d} \frac{ a_{k} }{a_{k-d}}
=\frac{2^{d+1}\pi^{d/2}\Gamma\left(\frac{k+1}{2}\right)
}{d!\Gamma\left(\frac{k+1-d}{2}\right)}
\end{equation}
if $k-1-d$ is even, and zero otherwise, $d=0,\dots,k-1$.

\bibliographystyle{ims}
\bibliography{bib/cone,bib/General}

\begin{thebibliography}{45}
\expandafter\ifx\csname natexlab\endcsname\relax\def\natexlab#1{#1}\fi
\expandafter\ifx\csname url\endcsname\relax
  \def\url#1{\texttt{#1}}\fi
\expandafter\ifx\csname urlprefix\endcsname\relax\def\urlprefix{URL }\fi
\providecommand{\eprint}[2][]{\url{#2}}

\bibitem[{Adler(1981)}]{Adler:1981}
\textsc{Adler, R.~J.} (1981).
\newblock \textit{The Geometry of Random Fields}.
\newblock John Wiley \& Sons, Chichester.

\bibitem[{Adler(2000)}]{Adler:2000}
\textsc{Adler, R.~J.} (2000).
\newblock On excursion sets, tube formulae, and maxima of random fields.
\newblock \textit{Annals of Applied Probability}, \textbf{10} 1--74.

\bibitem[{Adler and Taylor(2007)}]{RFG}
\textsc{Adler, R.~J.} and \textsc{Taylor, J.~E.} (2007).
\newblock \textit{Random fields and their geometry}.
\newblock Birkh\"{a}user, Boston.

\bibitem[{Becker et~al.(2009)Becker, Bobin and Cand\`es}]{nesta}
\textsc{Becker, S.}, \textsc{Bobin, J.} and \textsc{Cand\`es, E.~J.} (2009).
\newblock {NESTA}: a fast and accurate first-order method for sparse recovery.
\newblock \textit{SIAM J. on Imaging Sciences}, \textbf{4} 1--39.

\bibitem[{Birnbaum(1954)}]{Birnbaum:1954}
\textsc{Birnbaum, A.} (1954).
\newblock Combining independent tests of significance.
\newblock \textit{Journal of the American Statistical Society}, \textbf{49}
  559--574.

\bibitem[{Calhoun et~al.(2004)Calhoun, Stevens, Pearlson and
  Kiehl}]{Calhoun:2004}
\textsc{Calhoun, V.}, \textsc{Stevens, M.}, \textsc{Pearlson, G.} and
  \textsc{Kiehl, K.} (2004).
\newblock {fMRI} analysis with the general linear model: removal of
  latency-induced amplitude bias by incorporation of hemodynamic derivative
  terms.
\newblock \textit{Neuro{I}mage}, \textbf{22} 252--257.

\bibitem[{Cao and Worsley(1999{\natexlab{a}})}]{Cao:Worsley:1999:Correlation}
\textsc{Cao, J.} and \textsc{Worsley, K.} (1999{\natexlab{a}}).
\newblock The geometry of correlation fields with an application to functional
  connectivity of the brain.
\newblock \textit{Annals of Applied Probability}, \textbf{9} 1021--1057.

\bibitem[{Cao and Worsley(1999{\natexlab{b}})}]{Cao:Worsley:1999:Hotelling}
\textsc{Cao, J.} and \textsc{Worsley, K.~J.} (1999{\natexlab{b}}).
\newblock The detection of local shape changes via the geometry of
  {H}otelling's ${T}\sp 2$ fields.
\newblock \textit{Annals of Statistics}, \textbf{27} 925--942.

\bibitem[{Carbonell and Worsley(2007)}]{Carbonell:Wilks}
\textsc{Carbonell, F.} and \textsc{Worsley, K.} (2007).
\newblock The geometry of the {W}ilks's $\lambda$ random field.
\newblock \textit{Annals of the institute of Statistical Mathematics}.
\newblock Submitted.

\bibitem[{Cohen and Sackrowitz(1993)}]{Cohen:Sack:Inadmiss}
\textsc{Cohen, A.} and \textsc{Sackrowitz, H.~B.} (1993).
\newblock Inadmissibility of studentized tests for normal order restricted
  models.
\newblock \textit{Annals of Statistics}, \textbf{21} 746--752.

\bibitem[{Friedman et~al.(2007)Friedman, Hastie, Hofling and
  Tibshirani}]{friedman_pathwise_2007}
\textsc{Friedman, J.~H.}, \textsc{Hastie, T.}, \textsc{Hofling, H.} and
  \textsc{Tibshirani, R.} (2007).
\newblock Pathwise coordinate optimization.
\newblock \textit{Annals of Applied Statistics}, \textbf{1} 302--332.

\bibitem[{Friman et~al.(2003)Friman, Borga, Lundberg and
  Knutsson}]{Friman:2003}
\textsc{Friman, O.}, \textsc{Borga, M.}, \textsc{Lundberg, P.} and
  \textsc{Knutsson, H.} (2003).
\newblock Adaptive analysis of {fMRI} data.
\newblock \textit{Neuro{I}mage}, \textbf{19} 837--845.

\bibitem[{Friston et~al.(1998)Friston, Fletcher, Josephs, Holmes, Rugg and
  Turner}]{Friston:1998}
\textsc{Friston, K.}, \textsc{Fletcher, P.}, \textsc{Josephs, O.},
  \textsc{Holmes, A.}, \textsc{Rugg, M.} and \textsc{Turner, R.} (1998).
\newblock Event-related {fMRI}: {C}haracterising differential responses.
\newblock \textit{Neuro{I}mage}, \textbf{7} 30--40.

\bibitem[{Friston et~al.(1995)Friston, Holmes, Worsley, Poline, Fritn and
  Frackowiak}]{friston_statistical_1995}
\textsc{Friston, K.~J.}, \textsc{Holmes, A.~P.}, \textsc{Worsley, K.~J.},
  \textsc{Poline, J.~P.}, \textsc{Fritn, C.~D.} and \textsc{Frackowiak, R.~S.}
  (1995).
\newblock Statistical parametric maps in functional imaging a general linear
  approach.
\newblock \textit{Human Brain Mapping}, \textbf{2} 189--210.

\bibitem[{Johansen and Johnstone(1990)}]{Johansen:Johnstone:1990}
\textsc{Johansen, S.} and \textsc{Johnstone, I.} (1990).
\newblock Hotelling's theorem on the volume of tubes: some illustrations in
  simultaneous inference and data analysis.
\newblock \textit{Annals of Statistics}, \textbf{18} 652--684.

\bibitem[{Johnstone and Siegmund(1989)}]{johnstone_hotellings_1989}
\textsc{Johnstone, I.} and \textsc{Siegmund, D.} (1989).
\newblock On hotelling's formula for the volume of tubes and naiman's
  inequality.
\newblock \textit{The Annals of Statistics}, \textbf{17} 184–194.

\bibitem[{Knowles and Siegmund(1989)}]{Siegmund:Knowles:1989}
\textsc{Knowles, M.} and \textsc{Siegmund, D.} (1989).
\newblock On {H}otelling's approach to testing for a nonlinear parameter in a
  regression.
\newblock \textit{International Statistical Review}, \textbf{57} 205--220.

\bibitem[{Lawson and Hanson(1995)}]{Lawson:Hanson}
\textsc{Lawson, C.~L.} and \textsc{Hanson, R.~J.} (1995).
\newblock \textit{Solving Least Squares Problems}.
\newblock Society for Industrial and Applied Mathematics, Philadelphia.

\bibitem[{Lin and Lindsay(1997)}]{Lin:Lindsay:1997}
\textsc{Lin, Y.} and \textsc{Lindsay, B.~G.} (1997).
\newblock Projections on cones, chi-bar squared distributions, and {W}eyl's
  formula.
\newblock \textit{Statistics \& Probability Letters}, \textbf{32} 367--376.

\bibitem[{Nardi et~al.(2008)Nardi, Siegmund and Yakir}]{yakir}
\textsc{Nardi, Y.}, \textsc{Siegmund, D.~O.} and \textsc{Yakir, B.} (2008).
\newblock The distribution of maxima of approximately {G}aussian random fields.
\newblock \textit{Annals of Statistics}, \textbf{36}.

\bibitem[{Perlman and Wu(1999)}]{Emperor:New:Tests}
\textsc{Perlman, M.~D.} and \textsc{Wu, L.} (1999).
\newblock The {E}mperor's new tests.
\newblock \textit{Statistical Science}, \textbf{14} 355--381.

\bibitem[{Pilla(2006)}]{Pilla:2006}
\textsc{Pilla, R.~S.} (2006).
\newblock Inference under convex cone alternatives for correlated data.
\newblock \textit{E-print}.
\newblock ArXiv:math/0506522v3.

\bibitem[{Polzehl and Tabelow(2006)}]{fmriR}
\textsc{Polzehl, J.} and \textsc{Tabelow, K.} (2006).
\newblock Analysing f{MRI} experiments with the fmri package in {R}. {V}ersion
  1.0 - {A} users guide.
\newblock \textit{Weierstrass Institute for Applied Analysis and Stochastics
  Technical Report}, \textbf{10}.

\bibitem[{Robertson et~al.(1988)Robertson, Wright and Dykstra}]{Rob:Wri:Dyk}
\textsc{Robertson, T.}, \textsc{Wright, F.~T.} and \textsc{Dykstra, R.~L.}
  (1988).
\newblock \textit{Order Restricted Statistical Inference}.
\newblock Wiley, New York.

\bibitem[{Shafie et~al.(2003)Shafie, Sigal, Siegmund and Worsley}]{Shafie}
\textsc{Shafie, K.}, \textsc{Sigal, B.}, \textsc{Siegmund, D.~O.} and
  \textsc{Worsley, K.~J.} (2003).
\newblock Rotation space random fields with an application to f{MRI} data.
\newblock \textit{Annals of Statistics}, \textbf{31} 1732--1771.

\bibitem[{Siegmund and Worsley(1995)}]{Siegmund:Worsley:1995}
\textsc{Siegmund, D.~O.} and \textsc{Worsley, K.~J.} (1995).
\newblock Testing for a signal with unknown location and scale in a stationary
  {G}aussian random field.
\newblock \textit{Annals of Statistics}, \textbf{23} 608--639.

\bibitem[{Sun(1993)}]{Sun:1993}
\textsc{Sun, J.} (1993).
\newblock Tail probabilities of the maxima of {G}aussian random fields.
\newblock \textit{Annals of Probability}, \textbf{21} 34--71.

\bibitem[{Sun and Loader(1994)}]{Sun:Loader:1994}
\textsc{Sun, J.} and \textsc{Loader, C.~R.} (1994).
\newblock Simultaneous confidence bands for linear regression and smoothing.
\newblock \textit{Annals of Statistics}, \textbf{22} 1328--1345.

\bibitem[{Sun et~al.(2000)Sun, Loader and McCormick}]{Sun:2000}
\textsc{Sun, J.}, \textsc{Loader, C.~R.} and \textsc{McCormick, W.~P.} (2000).
\newblock Confidence bands in generalized linear models.
\newblock \textit{Annals of Statistics}, \textbf{28} 429--460.

\bibitem[{Takemura and Kuriki(1997)}]{Takemura:Kuriki:Chi:Bar}
\textsc{Takemura, A.} and \textsc{Kuriki, S.} (1997).
\newblock Weights of {$\overline{\chi}{}\sp 2$} distribution for smooth or
  piecewise smooth cone alternatives.
\newblock \textit{Annals of Statistics}, \textbf{25} 2368--2387.

\bibitem[{Takemura and Kuriki(2002)}]{Takemura:Kuriki:Equivalence:2002}
\textsc{Takemura, A.} and \textsc{Kuriki, S.} (2002).
\newblock On the equivalence of the tube and {E}uler characteristic methods for
  the distribution of the maximum of {G}aussian fields over piecewise smooth
  domains.
\newblock \textit{Annals of Applied Probability}, \textbf{12} 768--796.

\bibitem[{Taylor(2006)}]{Taylor:2002:Kinematic}
\textsc{Taylor, J.~E.} (2006).
\newblock A {G}aussian kinematic formula.
\newblock \textit{Annals of Probability}, \textbf{34} 122--158.

\bibitem[{Taylor and Adler(2003)}]{Taylor:Adler:2003:EC:Manifolds}
\textsc{Taylor, J.~E.} and \textsc{Adler, R.~J.} (2003).
\newblock Euler characteristics for {G}aussian fields on manifolds.
\newblock \textit{Annals of Probability}, \textbf{31} 533--563.

\bibitem[{Taylor et~al.(2005)Taylor, Takemura and Adler}]{TTA:Validity}
\textsc{Taylor, J.~E.}, \textsc{Takemura, A.} and \textsc{Adler, R.~J.} (2005).
\newblock Validity of the expected {E}uler characteristic heuristic.
\newblock \textit{Annals of Probability}, \textbf{33} 1362--1396.

\bibitem[{Taylor and Vadlamani(2011)}]{taylor_vadlamani}
\textsc{Taylor, J.~E.} and \textsc{Vadlamani, S.} (2011).
\newblock Random fields and the geometry of {W}iener space.
\newblock \textit{Annals of Probability}.
\newblock To appear., \urlprefix\url{http://arxiv.org/abs/1105.3839}.

\bibitem[{Taylor and Worsley(2007)}]{TW:noniso}
\textsc{Taylor, J.~E.} and \textsc{Worsley, K.~J.} (2007).
\newblock Detecting sparse signals in random fields, with an application to
  brain mapping.
\newblock \textit{Journal of the American Statistical Association},
  \textbf{102} 913--928.

\bibitem[{Taylor and Worsley(2008)}]{Taylor:Worsley:Roy}
\textsc{Taylor, J.~E.} and \textsc{Worsley, K.~J.} (2008).
\newblock Random fields of multivariate test statistics, with applications to
  shape analysis.
\newblock \textit{The Annals of Statistics}, \textbf{36} 1--27.

\bibitem[{Worsley et~al.(1999)Worsley, Andermann, Koulis, MacDonald and
  Evans}]{Worsley:Nonisotropic:1999}
\textsc{Worsley, K.}, \textsc{Andermann, M.}, \textsc{Koulis, T.},
  \textsc{MacDonald, D.} and \textsc{Evans, A.} (1999).
\newblock Detecting changes in nonisotropic images.
\newblock \textit{Human Brain Mapping}, \textbf{8} 98--101.

\bibitem[{Worsley et~al.(2002)Worsley, Liao, Aston, Petre, Duncan, Morales and
  Evans}]{Worsley:fmristat}
\textsc{Worsley, K.}, \textsc{Liao, C.}, \textsc{Aston, J.}, \textsc{Petre,
  V.}, \textsc{Duncan, G.}, \textsc{Morales, F.} and \textsc{Evans, A.} (2002).
\newblock A general statistical analysis for {fMRI} data.
\newblock \textit{Neuro{I}mage}, \textbf{15} 1--15.

\bibitem[{Worsley and Taylor(2006)}]{Worsley:Taylor:2005:latency}
\textsc{Worsley, K.} and \textsc{Taylor, J.} (2006).
\newblock Detecting {fMRI} activation allowing for unknown latency of the
  hemodynamic response.
\newblock \textit{Neuro{I}mage}, \textbf{29} 649--654.

\bibitem[{Worsley(1994)}]{Worsley:1994:Chi:t:F}
\textsc{Worsley, K.~J.} (1994).
\newblock Local maxima and the expected {E}uler characteristic of excursion
  sets of $\chi\sp 2,\ {F}$ and $t$ fields.
\newblock \textit{Advances in Applied Probability}, \textbf{26} 13--42.

\bibitem[{Worsley(1995{\natexlab{a}})}]{Worsley:1995:Boundary}
\textsc{Worsley, K.~J.} (1995{\natexlab{a}}).
\newblock Boundary corrections for the expected {E}uler characteristic of
  excursion sets of random fields, with an application to astrophysics.
\newblock \textit{Advances in Applied Probability}, \textbf{27} 943--959.

\bibitem[{Worsley(1995{\natexlab{b}})}]{worsley_estimating_1995}
\textsc{Worsley, K.~J.} (1995{\natexlab{b}}).
\newblock Estimating the number of peaks in a random field using the hadwiger
  characteristic of excursion sets, with applications to medical images.
\newblock \textit{The Annals of Statistics}, \textbf{23} 640–669.

\bibitem[{Worsley(2001)}]{Worsley:Scale:Chisq:2001}
\textsc{Worsley, K.~J.} (2001).
\newblock Testing for signals with unknown location and scale in a $\chi^2$
  random field, with an application to {fMRI}.
\newblock \textit{Advances in Applied Probability}, \textbf{33} 773--793.

\bibitem[{Worsley et~al.(1996)Worsley, Marrett, Neelin, Vandal, Friston and
  Evans}]{worsley_unified_1996-1}
\textsc{Worsley, K.~J.}, \textsc{Marrett, S.}, \textsc{Neelin, P.},
  \textsc{Vandal, A.~C.}, \textsc{Friston, K.~J.} and \textsc{Evans, A.~C.}
  (1996).
\newblock A unified statistical approach for determining significant signals in
  images of cerebral activation.
\newblock \textit{Human Brain Mapping}, \textbf{4} 58–73.

\end{thebibliography}

\end{document}